\RequirePackage{fix-cm}
\pdfoutput=1
\documentclass[leqno, centertags, 12pt]{article}

\usepackage{amsmath}
\usepackage{amscd}
\usepackage{amssymb}
\usepackage{verbatim}

\usepackage[T1]{fontenc}
\usepackage{ae} 
\usepackage{aecompl}
\usepackage{latexsym}
\usepackage{manfnt}
\usepackage[sans]{dsfont}
\usepackage{palatino}
\usepackage[sc, osf]{mathpazo} 
\usepackage{eulervm}

\usepackage{fullpage}

\usepackage{graphicx}
\newcommand{\minus}{\,{\scalebox{0.75}[1.0]{$-$}}\hspace{0.2em}}

\DeclareMathOperator{\genus}{{\mathfrak g}\hspace{0.5pt}}
\DeclareMathOperator{\Aut}{{Aut}\hspace{0.5pt}}

\DeclareMathOperator{\Mod}{{Mod}\hspace{0.5pt}}
\DeclareMathOperator{\T}{{\mathcal I}\hspace{0.5pt}}
\DeclareMathOperator{\Sp}{{Sp}\hspace{0.5pt}}

\DeclareMathOperator{\ttg}{{\mathfrak Tb}\hspace{0.5pt}} 
\DeclareMathOperator{\ccc}{{\mathcal C}\hspace{0.5pt}} 
\DeclareMathOperator{\link}{{\mathfrak L}\hspace{0.5pt}}

\DeclareMathOperator{\mme}{{EMod}\hspace{0.5pt}}
\DeclareMathOperator{\spe}{{ESp}\hspace{0.5pt}}
\DeclareMathOperator{\comm}{{Comm}\hspace{0.5pt}}

\newcommand{\Z}{{\bf Z}}

\newcommand{\rr}{{\bf R}}

\newcommand{\gen}{{\mathfrak g}}

\newcommand{\ffc}{F_\mathrm{cc}}
\newcommand{\fft}{F_\mathrm{to}}

\newcommand{\tto}{\,{\to}\,}

\newcommand{\tts}{\ttg(S)}
\newcommand{\ccs}{ \ccc(S)}
\newcommand{\tors}{\T(S)}

\newcommand{\seq}[4]{{#1}_1 \,{=}\,{#2}\/\/,\, {#1}_2\/\/,\,\ldots\,, {#1}_{#4}\,{=}\,{#3}}

\newcommand{\sep}{_{\mathrm{sep}}}

\newcommand{\langles}{\mathbin{\langle}}  
\newcommand{\rangles}{\mathbin{\rangle}}  
\newcommand{\vv}[1]{\langles #1 \rangles} 
\newcommand{\llvv}[1]{\langles\hspace{0.05em} #1 \hspace{0.07em}\rangles} 
\newcommand{\aavv}[1]{\langle\hspace{0.07em} {#1}\hspace{0.05em}\rangle} 
\newcommand{\geqs}{\mathbin{\geq}} 
\newcommand{\leqs}{\mathbin{\leq}}

\newcommand{\setmiss}{{\setminus}\hspace{0.05em}}
\newcommand{\llsetmiss}{\,\setmiss}

\newcommand{\elem}{\,{\in}\,}

\newcommand{\eeq}{{\hskip.05em\,{=}\,}}

\def\sss{\hskip.05em\ }

\def\dds{\hskip.1em\ }
\def\trs{\hskip.15em\ }

\newcommand{\delete}[1]{}

\makeatletter
\renewcommand{\@makefntext}[1]{\parindent=0em\noindent
\hbox to 0.4em{\hss\@makefnmark}#1}
\makeatother

\parskip=\the\bigskipamount

\newcommand{\proof}{\emph{Proof.}\hspace{0.3em}}
\newcommand{\eproof}{ $\blacksquare$}
\newcommand{\esubproof}{ $\square$}

\newcounter{mysectionnumber}
\newcommand{\mysection}[2]{\setcounter{footnote}{0}\setcounter{equation}{0}\refstepcounter{mysectionnumber}
\vspace{-8.5ex}
\section*{\begin{center}\textnormal{{\themysectionnumber.} {#1}}\end{center}}\label{#2}\vspace{0ex}}

\newcounter{myparnum}[mysectionnumber]
\newcommand{\mypar}[2]{\refstepcounter{myparnum}{\vspace{-\bigskipamount} \paragraph{\textit{{\themyparnum. #1}\label{#2}}} \hspace{-0.7em}}}
\renewcommand{\themyparnum}{\themysectionnumber.\arabic{myparnum}}

\newcounter{mylemmanum}[myparnum]
\renewcommand{\themylemmanum}{\themyparnum.\arabic{mylemmanum}}
\newcommand{\mylemma}[2]{\refstepcounter{mylemmanum}{\vspace{-\bigskipamount}\paragraph{\textit{\themylemmanum. #1}\label{#2}}\hspace{-0.7em}}}

\numberwithin{equation}{section}

\newcommand{\myitpar}[1]{\vspace{-\bigskipamount}\paragraph{\textit{#1}}\hspace{-0.7em}}

\newcommand{\mytitle}[1]{\textbf{\textit{#1}}}

\newcommand{\mynonumbersection}[2]{
\vspace{-0.0ex}
\section*{{}\hspace*{0.00em}$\phantom{1.}$\textnormal{{#1}}}\label{#2}\vspace{\bigskipamount}\vspace{-0.0ex}}  

\title{Torelli buildings and their automorphisms\vspace*{0ex}}
\author{Benson Farb and Nikolai V. Ivanov}
\date{}

\begin{document}

\maketitle

\footnotetext{\hspace{-0.45em}The work reported in this paper was largely completed about ten years ago.\vspace{\medskipamount}}

\footnotetext{\hspace{-0.45em}The work of Benson Farb was supported in part by NSF grants DMS-9704640 
and DMS-0244542.} 
\footnotetext{\hspace{-0.45em}The work of Nikolai V. Ivanov was supported in part by the University of Chicago.}

\vspace{9ex}

\mynonumbersection{Preface}{preface}

It is well known that each compact orientable surface $R$ gives rise to a simplicial complex $\ccc(R)$
playing in the theory of Teichm{\"u}ller modular groups a role similar to the role of Tits buildings in theories of algebraic and arithmetic groups.
In this paper we introduce, for each closed orientable surface $S$, 
an analogue of Tits buildings more suitable for investigation of the Torelli group $\tors$ of $S$.
It is a simplicial complex with some additional structure.
We denote this complex with its additional structure by $\tts$ and call it the {\em Torelli building\/} of $S$. 
The main result of this paper shows that Torelli buildings of surfaces of genus $\geqs 5$ have only obvious automorphisms, 
and identifies the group of automorphisms of $\tts$. 
Namely, we prove that for such surfaces $S$ every automorphism of $\tts$ is induced by a diffeomorphism of $S$.
This theorem is an analogue of the theorem of the second author \cite{i-ihes,i-imrn} about authomorphisms of complexes of curves $\ccc(R)$.
The proof is based on this theorem about automorphisms of complexes of curves, and uses methods inspired by its proof \cite{i-ihes,i-imrn}.

Our theorem about automorphisms of Torelli buildings is intended for applications to automorphisms and virtual automorphisms of Torelli groups.
These applications are also inspired by the results and methods of \cite{i-ihes,i-imrn}, and are proved following the scheme developed
by the second author \cite{i0,i-aut,i-ihes} for studying automorphisms of Teichm{\"u}ller modular groups.
Compared to Teichm{\"u}ller modular groups, Torelli groups are much less accessible and much less is known about them.
 
Our main result about Torelli groups is the identification of all isomorphisms between subgroups of finite index in the Torelli group $\tors$ of a closed surface $S$. Namely, every such isomorphism is equal to (the restriction of) to a conjugation by an element of the extended mapping class group $\mme (S)$ of $S$. 
This result immediately allows to find the automorphism group and the outer automorphism group of $\tors$, and to find its abstract commensurator. 
The deduction of our theorems about Torelli groups from the theorem about automorphisms of Torelli building 
is independent of the methods of the present paper and uses only its main result.
It will be published on some other occasion.

The results of this paper, together with their applications to Torelli groups, were announced in \cite{fi-announ}.

\mysection{Introduction}{intro}

The notations introduced in this section will be used throughout the whole paper.
For a compact connected surface $R$ we will denote by $\genus(R)$ its genus.
Throughout the paper we will denote by $S$ a closed orientable surface and by $\gen\eeq\genus(S)$ its genus.

The main resuts of this paper are proved under the assumption $\gen\eeq\genus(S)\geqs 5$. 
But most of our proofs are valid under weaker assumptions $\gen\geqs 4$ or  $\gen\geqs 3$. 
It is quite reasonable to believe that the assumption $\gen\geqs 5$ can be weakened
to $\gen\geqs 4$ throughout the whole paper at the cost of additional technicalities.
The authors were not really interested in proving the main result with the
weakest possible restrictions on $\gen$, and preferred the most natural (at least for
the authors) arguments. The latter always worked for $\gen\geqs 5$.

The \emph{Teichm{\"u}ller modular group\footnote{Admittedly, another term, uniformative and colorless, 
is more widely used.} $\Mod(S)$ of\/ $S$} is defined
as the group of isotopy classes of orientation-preserving diffeomorphisms $S\tto S$.
We will need also the \emph{extended Teichm{\"u}ller modular group\/} $\mme(S)$
of $S$ which is defined as the group of the isotopy classes of \emph{all\/} diffeomorphisms
$S\tto S$. Let us fix an orientation of $S$. Then
the algebraic intersection number provides a nondegenerate,
skew-symmetric, bilinear form on $H\eeq H_1(S,\Z)$, called the {\em
intersection form.} The natural action of $\Mod(S)$ on $H$ preserves 
the intersection form. If fix a symplectic basis in $H$, then
we can identify the group of symplectic automorphisms of $H$
with the integral symplectic group $\Sp(2\gen,\Z)$ and the
action of $\Mod(S)$ on $H$ leads to a natural surjective homomorphism $\Mod(S)\tto \Sp(2\gen,\Z)$.
 
The {\em Torelli group} of $S$, denoted by $\tors$, is defined as the kernel of this homorphism; 
that is, $\tors$ is the subgroup of $\Mod(S)$ consisting of elements which act trivially on $H_1(S,\Z)$. 
The groups $\tors$ plays an important role in algebraic geometry and low-dimensional topology. 
At the same time $\tors$ is an interesting group in its own right. 

The modern theory of Torelli groups starts with the work of D. Johnson in the 1980-ies 
(see, for example, surveys \cite{Jo1} and \cite{Ha1}).
The first spectacular results of D. Johnson (preceded by some preliminary technical results)
is his theorem \cite{Jo2} to the effect that the group 
$\tors$ is finitely generated if the genus $\genus(S)\geqs 3$. 
Thirty years later it is not known whether or not $\tors$ admits a finite presentation for $\genus(S)\geqs 3$. 
It seems that D. Johnson's work remains to be poorly understood: we understand his proofs, but we do not understand why they work.
Unfortunately for mathematicians interested in Torelli groups, but reportedly not for himself, D. Johnson left mathematics in the late 1980-ies.
He left the academic world already in the late 1970-ies.

A natural problem is to find all automorphisms of groups $\tors$ or, better, all isomorphisms between subgroups of finite index of $\tors$.
This problem motivated this paper, which devoted to the first part of a solution.
As a motivation for the reader, let us outline our main results in this direction.

As usual, for any (discrete) group $G$ we denote by $\Aut(G)$ the group of all automorphisms of $G$.
If $H$ is a normal subgroup of a group $G$, then $G$ acts on $H$ by conjugations.
Namely, if $g\elem G$, then the formula $h\mapsto g^{-1} h g$, where $h\elem H$, defines and automorphism $c_g\colon H\tto H$. 
One immediately checks that the map $g\mapsto c_g$ is a group homomorphism, which we will denote by $\mbox{co}\colon G\tto \Aut(G)$.

Since $\tors$ is a normal subgroup of $\mme(S)$, the Teichm{\"u}ller modular group $\mme(S)$ acts on $\tors$ by conjugation. 
In particular, the homomorphism $\mbox{co}\colon\mme(S)\tto \Aut(\tors)$ is defined. 
The main result about automorphisms is the following theorem.
 
\mypar{Theorem.}{aut-group} \emph{If\sss $\gen\eeq\genus(S)\geqs 5$, then\sss  $\textup{\mbox{co}}\colon\mme(S)\tto \Aut(\tors)$\sss is an isomorphism.}

Note that the conclusion of the Theorem \ref{aut-group} is false if $\gen\eeq 2$, 
since in this case $\tors$ is a countably generated free group by a theorem of Mess \cite{Me}. 
Theorem \ref{aut-group} and its proof are inspired by the work of 
the second author \cite{i0,i-aut} about automorphisms of\sss $\Mod(S)$, and J. McCarthy's \cite{mc-aut} version of it.
The methods of papers \cite{i0,i-aut}, \cite{mc-aut} are not powerful enough to deal with automorphisms of Torelli groups.
In our work we followed the methods developed later by the second author \cite{i-ihes,i-imrn}. 
These methods automatically lead to stronger results and allow to deal not only with automorphisms of $\tors$, 
but also with isomorphisms between subgroups of finite index of $\tors$.
Namely, they lead to the following theorem.

\mypar{Theorem.}{virtual-aut} \emph{If\sss $\gen\eeq\genus(S)\geqs 5$, then every isomorphism\sss $\Gamma_1 \tto \Gamma_2$\sss 
between two subgroups\sss $\Gamma_1$,\sss $\Gamma_2$\sss of finite index in\sss\ $\tors$\sss is induced by the conjugation\sss $c_g$\sss by an element\dds $g\elem\mme(S)$.}

This theorem allows immediately find the so-called \emph{abstract commensurator} of $\tors$.
In the present context it is most natural to refer the reader to \cite{i-imrn} for a definition. 

\mypar{Theorem.}{abs-comm} \emph{If\dds  $\gen\eeq\genus(S)\geqs 5$, then the abstract commensurator\trs 
$\comm (\tors)$\sss of\trs $\tors$\sss is isomorphic to\sss $\spe(2\gen, \Z)$,\sss
where\sss $\spe (2\gen, \Z)$\sss is the group of automorphims of\sss $\Z^{2g}$\sss preserving the standard symplectic form up to a sign.}

\mytitle{Remark.} The \emph{extended symplectic group}\dds $\spe (2\gen, \Z)$ contains the standard symplectic group $\Sp(2\gen,\Z)$ 
as a subgroup of index $2$.

An easy application of this result is a new proof of the following special case of a theorem of the second author \cite{i-ar}.

\mypar{Theorem.}{non-arithmetic} \emph{If\/  $\gen\eeq\genus(S)\geqs 5$, then 
the Torelli group\/ $\tors$ is not isomorphic to any arithmetic group.}

While this theorem is neither new, nor the most general one of this type,
our new proof shed an additional light on it (cf. \cite {i-imrn}, Theorem 6 and the discussion following it).

A key ingredient in the proof of Theorems \ref{aut-group} and \ref{virtual-aut} is an analogue of the theorem about automorphisms
of the complexes of curves due to the second author \cite{i-ihes,i-imrn}.
Recall that the \emph{complex of curves}\/ $ \ccc(R)$ of an orientable surface $R$ is the simplicial complex having as the set of vertices the set of isotopy
classes of all non-trivial (i.e. not bounding a disc in $R$ and not homotopic to a component of the boundary $\partial R$) circles in $R$;
a set of vertices forms a simplex if they can be represented by disjoint circles. 
See the subsections \ref{s-complexes} and \ref{cc} for more details about simplicial complexes and complexes of curves,
including the definition of simplicial complexes.

An analogue of the complex of curves in the theory of Torelli groups is a simplicial complex having as vertices the isotopy
classes of separating circles on $R$ and pairs of the isotopy classes of non-separating circles $C$, $C'$ on $R$ 
such that the union $C\cup C'$ is separating. 
It is convenient to equip this simplicial complex equipped with some additional structure.
See the subsection \ref{tg-def} for the definition. 
This simplicial complex together with its additional structure is called the \emph{Torelli building of}\/ $R$ and is denoted by $ \ttg(R)$.

\vspace{-2.3ex}
\paragraph{\textit{Main Theorem.}} \emph{If\dds $\gen{\eeq}\genus(S)\geqs 5$,\sss then every automorphism of the 
Torelli building\sss $\tts$\sss of\sss $S$\sss is induced by a diffeomorphism of\sss $S$.}
 
The Main Theorem is an analogue of the theorem of the second author \cite{i-ihes,i-imrn} about automorphisms of the complexes of curves.
The  the definition of the Torelli building $\tts$ will be given in Section \ref{t-geom},
and after this Main Theorem will be stated again in as Theorem \ref{automorphisms of geometry}.
Sections \ref{proof-conn-ths}--\ref{en-edges-aut} are devoted to its proof. 
Section \ref{prelim} is devoted to an overview of basic definitions and notations used in this paper.  
\vspace{1ex}

\mysection{Preliminaries}{prelim}

This section is devoted to basic definitions and facts used throughout the paper.

\mypar{Surfaces.}{surfaces} By a \emph{surface} we understand a compact orientable $2$-manifold with boundary.
The boundary is allowed to be empty, and surfaces with empty boundary are called \emph{closed}\/ surfaces. 
Our main results are concerned only with closed surfaces, 
but surfaces with non-empty boundary inevitably show up in the proofs. 

By a \emph{subsurface}\/ of a surface $R$ we understand a codimension $0$ submanifold
$Q$ of $R$ such that each component of the boundary $\partial Q$ is either equal to
a component of $\partial R$, or disjoint from $\partial R$. 
Obviously, if $R$ is a closed connected surface and $Q$ is a subsurface of $R$,
then either $Q {\eeq} R$, or the boundary $\partial R$ is non-empty.

\mypar{Circles.}{circles}
The most basic tool in the topology of surfaces are submanifolds of surfaces 
diffeomorphic to the standard circle $S^1$ in $\rr^2$. These submanifolds are
assumed to be disjoint from the boundary (if the boundary is non-empty).
We we will call such a submanifold of a surface $R$ a \emph{circle in\/ $R$}.
A circle in $R$ is called \emph{non-peripheral}\/ if does not bound 
an annulus together with a component of the boundary $\partial R$. 
If $R$ is a closed surface, then, obviously, all circles in $R$ are non-peripheral. 
A circle called \emph{non-trivial}\/ if it is \emph{non-peripheral}\/ and does not bound a disc in $R$. 
Alternatively, a \emph{non-peripheral}\/ circle can be defined as circle which 
is not homotopic to a component of the boundary $\partial R$.
Similarly, a \emph{non-trivial}\/ circle can be defined 
for closed surfaces as a cirle which cannot be deformed into a point, 
and for surfaces with non-empty boundary as a circle which cannot be deformed into the boundary.

\mypar{Separating circles.}{sep-circles} A circle $D$ in a connected surface $R$ 
is called \emph{separating}, if the result of cutting $R$ along $D$ is not connected, 
or, what is the same, $R\llsetmiss D$ is not connected. 
Then both the result of cutting and $R\llsetmiss D$ consist of two components.
One more way to define \emph{separating}\/ circles in $R$ is to define them as
\emph{homologically trivial}\/ circles, i.e. as circles in $R$ such that their 
fundamental class, considered as a $1$-dimensional homology class in $R$, is equal to $0$.
The last approach has the advantage of working in the expected way also for
non-connected surfaces. Separating circles are often called also \emph{bounding circles}.

A separating circle $D$ in a connected surface $R$ divides $R$ into two subsurfaces of $R$. 
These subsurfaces are equal to the closures of two components of $R\llsetmiss D$, 
and $D$ is a boundary component of each of them. If $R$ is closed, then $D$
is the boundary of both these subsurfaces. We will often call these subsurfaces 
\emph{the parts into which\/ $D$ divides $S$}. 

If $h$ is the genus of either of two subsurface bounded by $D$, 
then we will say that \emph{$D$ is a separating circle of genus\/ $h$}. 
While in general there are two possible values of $h$, this expression is unambiguos.

\mypar{Bounding pairs of circles.}{bp} A \emph{bounding pair of circles}\/ 
on a connected surface $R$ is defined as an unordered pair $C$, $C'$ of disjoint non-isotopic circles in $R$
such that both circles $C$, $C'$ are non-separating, but the result of cutting $R$
along $C\cup C'$ is not connected. Then it consists, obviously, of two components.
Alternatively, $C$, $C'$ is a \emph{bounding pair of circles}\/ if $R\llsetmiss (C\cup C')$ is not connected.

We should stress that a bounding pair of circles is a pair of submanifolds $C$, $C'$  of $R$, and
not the submanifold $C\cup C'$. Strictly speaking, a bounding pair of circles $C$, $C'$
is the $2$-element set $\{C,C'\}$. We will use the notation $\{C,C'\}$ only when this is essential,
and will use the term \emph{a bounding pair of circles $C$, $C'$} when this will cause no harm.
Of course, the pair $\{C,C'\}$ and the submanifold $C\cup C'$ carry the same information. 
But pairs have some (admittedly, minor) technical advantage.
See, for example, the subsection \ref{tg-def} and the proof of Lemma \ref{singling out} below. 

As in the case of separating circles, one can give also a homological definition.
It takes an especially simple form for closed surfaces. Namely, 
an unordered pair $C$, $C'$ of disjoint circles on a closed connected surface $R$ is a
\emph{bounding pair of circles}\/ if the fundamental class of the submanifold $C\cup C'$ with
some orientation, considered as a $1$-dimensional homology class in $R$, is equal to zero.
Note that $C\cup C'$ consists of two orientable components and hence admits $4$ orientations.
If $C$, $C'$ is a bounding pair of circles on a closed surface $R$, then for $2$ of these $4$ orientations
the fundamental class is equal to $0$, and for $2$ other is not. If $C$, $C'$ is
not a bounding pair of circles, then all $4$ orientations lead to the non-zero fundamental class.

The reader may skip the following paragraph and go directly to the subsection \ref{s-complexes}.

In order to give a homological definition of bounding pairs of circles in the general case 
of a possibly non-closed and non-connected surface $R$ it is convenient to fix an orientation 
of $R$ (recall that all our surfaces are orientable). This orientation induces an orientation
of the boundary $\partial R$, and hence an orientation of every boundary component of $R$.
By the fundamental class of a boundary component we will understand the fundamental class
of this component with this induced orientation. Now we can define a \emph{bounding pair of circles}\/ in $R$
as an unordered pair $C$, $C'$ of disjoint circles in $R$, such that
the fundamental class of $C\cup C'$ with \emph{some orientation}, 
considered as a $1$-dimensional homology class in $R$, is equal to the sum of the fundamental classes
of several (different) boundary components. Because we fixed an orientation of $R$, if
$C$, $C'$ is a bounding pair of circles, then only $1$ of the $4$ orientations of $C\cup C'$ satisfies
this condition. If we replace our fixed orientation of $R$ by the opposite orientation,
then orientations of all boundary components of $R$ will be replaced by the opposite orientations,
and then another orientation of $C\cup C'$ will satisfy our condition (obviously, the
orientation opposite to the first one). And if $C$, $C'$ is not a bounding pair of circles, then none
of the $4$ orientations of $C\cup C'$ will satisfy this condition.

\mypar{Isotopy classes of circles.}{isotopy-classes} By a standard abuse of language, 
the term \emph{circle}\/ is often used also for an \emph{isotopy class of circles} 
in the class of submanifolds of the surface in question. 
In this paper, the distinction between circles on surfaces and their isotopy classes is usually essential,
and we will avoid this abuse of notations to the extent possible. 

The isotopy classes of circles are vertices of the complex of curves, and, in the case
of separating circles in $S$, also of $\tts$. In view of this, these isotopy classes will
be called either \emph{vertices}\/ (of a specified complex), or just the isotopy classes.
For a circle $C$ in a surface $R$ we will denote \emph{the isotopy class of\/ $C$ in\/ $R$} by $\llvv{C}$.
The surface $R$ is usually clear from the context, even if $C$ is also a circle in some
other relevant surfaces (for example, some subsurfaces of $R$).

\mypar{Simplicial complexes.}{s-complexes} By a simplicial complex $\mathcal V$ we understand a pair consisting of 
a set $\mbox{V}$ together with a collection of \emph{finite subsets of}\/ $\mbox{V}$. 
As usual, we think of $\mathcal V$ as a structure on the set $\mbox{V}$,
namely a structure of a simplicial complex.
Elements of $\mbox{V}$ are called \emph{the vertices of\/ $\mathcal V$}, and subsets of $\mbox{V}$ 
from the given collection are called \emph{the simplices of\/ $\mathcal V$}.
These data are required to satisfy only one condition: \emph{a subset of a simplex is also a simplex}.
For every subset $\mbox{W}\subset \mbox{V}$, the structure $\mathcal V$ of a simplicial complex on $\mbox{V}$ leads to 
a canonical structure of a simplicial complex $\mathcal W$ on $\mbox{W}$. 
Namely, a subset of $\mbox{W}$ is declared to be a simplex of $\mathcal W$ 
if and only if it is a simplex of $\mathcal V$. The complex $\mathcal W$
is called the \emph{full subcomplex of\/ $\mathcal V$}\/ spanned by the set $\mbox{W}$.

For two simplicial comlexes $\mathcal{V}_1$, $\mathcal{V}_2$ with 
the sets of vertices $\mbox{V}_1$, $\mbox{V}_2$ respectively, 
an \emph{isomorphism}\/ $F\colon\mathcal{V}_1{\tto}\mathcal{V}_2$ is defined 
as a bijective map $F\colon\mbox{V}_1{\tto}\mbox{V}_2$ such that 
$\Delta$ is a simplex of $\mathcal{V}_1$ if and only if $F(\Delta)$ is a simplex of $\mathcal{V}_2$. 
Alternatively, we may requite that the images of simplices under both $F$ and $F^{-1}$ are also simplices.
As usual, we do not distinguish between an isomorphism $\mathcal{V}_1{\tto}\mathcal{V}_2$ 
and its underlying bijection $F\colon\mbox{V}_1{\tto}\mbox{V}_2$. 
The notion of a \emph{simplicial map}\/  $\mathcal{V}_1{\tto}\mathcal{V}_2$ 
which is not necessarily an isomorphism is well known, 
but is not needed for the purposes of this paper.

The \emph{dimension}\/ of a simplex is defined as its cardinality minus $1$.
An \emph{edge}\/ is defined as a simplex of dimension $1$, i.e. as a pair of vertices which is a simplex.
A \emph{triangle}\/ is defined as a simplex of dimension $2$.

We will say that two vertices $v$, $w$ are \emph{connected by an edge}\/
if either $\{v,w\}$ is a simplex of dimension $1$, or $v{\eeq}w$. In the second case $\{v,w\}{\eeq}\{v\}$ and
$\{v,w\}$ is a simplex of dimension zero.
We will say that $v$, $w$ are \emph{connected by a sequences of edges}\/ if there is a sequence
\[
\seq{v}{v}{w}{n}
\]
of vertices such that $v_i$ and $v_{i+1}$ are connected by an edge for every $i{\eeq}1,2,\ldots,n{\minus}1$.
A simplicial complex is called \emph{connected}\/ if every two vertices are connected by a sequence
of edges. Obviously, the property of being connected by a sequence of edges is an equivalence relation.
An equivalence class of vertices together with its induced structure of a simplicial complex
is called a \emph{component of connectedness}\/ of the original complex.

A simplicial complex is called a \emph{flag complex}\/ if the following condition holds:
\emph{a finite set $\Delta$ of vertices is a simplex if
and only if every subset of $\Delta$ with cardinality $2$ is a simplex}. A flag complex is determined by
its sets of vertices and edges, i.e. it is essentially a graph without loops and multiple edges.
All simplicial complexes considered in the paper are, in fact, flag complexes.

The higher dimensional simplices of a simplicial complex are 
exactly that part of the structure a simplicial complex
which makes simplicial complexes interesting and useful, even for flag complexes. 
In this paper we will work only with the vetices, edges, and triangles. 
But eventually we will have to refer to results which are proved 
by using higher dimensional simplices in an essential way. 
The most important of them is a theorem of the second 
author \cite{i-ihes,i-imrn} about automorphisms of complexes $\ccs$.

\mypar{The complex of curves.}{cc} Let $R$ be a compact surface, possibly with non-empty boundary.
The \emph{complex of curves}\/ $\ccc(R)$ is a simplicial complex in the sense of the subsection \ref{s-complexes}.
Its vertices are the isotopy classes $\llvv{C}$ of non-trivial circles $C$ in $R$. 
A collection of such isotopy classes is declared to be a simplex if and only if it is either empty, or
the isotopy classes from this collection can be represented by pair-wise disjoint circles. 
In other words, if $C_1$, $C_2$, \ldots, $C_n$ are pair-wise disjoint circles,
then the set $\{\vv{C_1}, \vv{C_2}, \ldots, \vv{C_n}\}$ is a simplex of $\ccc(R)$, 
and all simplices can we constructed in this way by the definition (the empty simplex corresponds to $n{\eeq}{\minus}1$).

If $C$ is a non-trivial separating (respectively, non-separating) circle on $R$,
we will call the corresponding vertex $\llvv{C}$ of $\ccc(R)$ a \emph{separating}\/ 
(respectively, \emph{non-separating}\/) vertex. 

\mytitle{$\ccc(R)$ is a flag complex}. This fact is well known to the experts.
For the benefit of the others, let us outline why this is the case. 

Let us restrict our attention by compact connected surfaces $R$ of negative Euler characteristic 
(the few exeptional surfaces with non-negative Euler characteristic can be easily dealt with directly).
It is well known that such a surface $R$ admits a hypebolic structure.
The latter can be defined as a riemannian metric on a surface having constant negative curvature
(which is usually required to be equal to ${\minus}1$, but this not important in our context).
Moreover, one can assume that the hyperbolic structure (riemannian metric) on $R$ has geodesic boundary.
Then every non-trivial circle in $R$ is isotopic to a unique geodesic circle in $R$,
and, moreover, if two circles are isotopic to some disjoint circles, then the geodesic circles
isotopic to them are either disjoint or equal.
Therefore, if we fix a hyperbolic structure on $R$, the geodesic circles may serve as canonical
representatives of isotopy classes of non-trivial circles, 
and a set $\{\gamma_1,\gamma_2,\ldots,\gamma_n\}$ of isotopy classes is a simplex of $\ccs$
if and only it the canonical representatives of $\gamma_1,\gamma_2,\ldots,\gamma_n$ 
are pair-wise disjoint. Clearly, this implies that $\ccc(R)$ is a flag complex.

\mytitle{Isomorphisms and automorphisms of complexes of curves.} Let
$F\colon\hspace{-1pt} R\tto Q$ be a diffeomorphism. 
Then $F$ induces an isomorphism $\ccc(R)\tto\ccc(Q)$ by the rule $\ffc\colon\aavv{C}\mapsto\aavv{F(C)}$. 
Clearly, $\ffc$ is a map from the set of vertices of $\ccc(R)$ to the set of vertices $\ccc(Q)$ which is bijective. 
Moreover, both $\ffc$ and $\ffc^{-1}$ take simplices to simplices.
We will denote the resulting automorphism by $\ffc\colon\ccc(R){\tto}\ccc(Q)$.
If $R{\eeq}Q$, then $\ffc$ is an automorphism of $\ccc(R)$. 
The following theorem of the second author \cite{i-ihes,i-imrn} tells us that there are no other automorphisms, at least if $\genus(R)\geq 2$.

\mypar{Theorem.}{aut-cc} \emph{If\/ $\genus(R)\geq 2$,
then every automorphism\/ $\ccc(R){\tto}\ccc(R)$ is equal to\/ $\ffc$\, for some diffeomorphism\/
$F\colon R{\tto} R$.}

This theorem was extended to the cases of $\genus(R){\eeq}1$ or $0$ by M. Korkmaz \cite{kork-th}, \cite{kork-pub}
and, independently but somewhat later, by F. Luo \cite{luo}. These extensions are not
needed for the present paper.

\mysection{Torelli buildings}{t-geom}

The \emph{Torelli building}\/ of $S$ is a simplicial complex with some
additional structure. By a standard abuse of language, we will denote by $\tts$
both this simplical complex and the Torelli building proper, i.e. the same complex
with this additional structure.

\mypar{The simplicial complex structure of $\tts$.}{tg-def} We start with 
the definitions the vertices and simplices of $\tts$.

\mytitle{Vertices.} There are two types of vertices of $\tts$.

The vertices of the first type are defined as the isotopy classes of non-trivial separating cirlces in $S$. 
These vertices are the vertices of the complex of curves $\ccs$ of the form $\langles D\rangles$,
where $D$ is a non-trivial separating circle in $S$. 
Following the conventions of the subsection \ref{cc} for the complexes of curves, 
we will call them \emph{separating vertices}. 
Occasionally we will call then $\mathcal{SC}$-\emph{vertices}.

The vertices of the second type are defined as the edges of $\ccs$ of the form
$\{\langles C_0\rangles, \langles C_1\rangles\}$, where $\{C_0,C_1\}$ is a bounding pairs of circles.
These vertices will be called \emph{bounding pairs},\/ or, occasionally, the $\mathcal{BP}$-\emph{vertices}.
Because the components of a bounding pair of circles are disjoint, 
and hence their isotopy classes are connected by an edge in $\ccs$,
every bounding pair of circles $P$ gives rise to a bounding pair 
consisting of the isotopy classes of its two elements. 
Extending our notations for the isotopy classes of circles, 
we will denote this bounding pair by $\llvv{P}$.\sss 
So, if $P{\eeq}\{C_0,C_1\}$,\sss then $\llvv{P}{\eeq}\{\vv{C_0},\vv{C_1}\}$.

\mytitle{Remark.} Our terminology does not allow using the term \emph{bounding pair}\/
as a short form of the term \emph{bounding pair of circles}.\/ This is a minor hassle for the authors ,
but, hopefully, not for the readers.

\mytitle{Simplices.} Given a set $\mbox{V}$ consisting of separating circles and
bounding pairs of circles in $S$, we will denote by $\llvv{\mbox{V}}$ the set $\{\llvv{A}\mid A\in \mbox{V}\hspace{0.1em}\}$.
In other terms, $\llvv{\mbox{V}}$ is the image of $\mbox{V}$ under the map $\bullet\mapsto\vv{\bullet}$.

Suppose that we are given a set $\mbox{W}_1$ of
separating circles in $S$ and a set $\mbox{W}_2$ of bounding pairs of circles in $S$.
Consider the set $\mbox{V}$ all circles in sight, 
i.e. the union of $\mbox{W}_1$ with 
the set of elements of bounding pairs of circles from $\mbox{W}_2$.
If all circles in the set $\mbox{V}$ are pair-wise disjoint, then the set
$\vv{\mbox{W}_1\cup \mbox{W}_2}$ is declared to be a \emph{simplex of}\/ $\tts$. 
Of course, sets $\mbox{W}_1$, $\mbox{W}_2$ are allowed to be empty. 
There are no other simplices in $\tts$.

\mytitle{${\tts}$ is a simplicial complex.} Clearly, a subset of a simplex is a simplex. 
Therefore, the set of vertices of both types together with the set of simplices forms a simplicial complex.
This is our simplicial complex $\tts$.

\mypar{The additional structure.}{add-structure} The additional structure 
on $\tts$ consists of two parts: the marking of each vertice by its type,
and a designation of some triangles as \emph{marked}\/ triangles. 
 
\mytitle{Markings of vertices.} The first part is a function 
from the set of vertices to the set $\{\mathcal{SC},\mathcal{BP}\}$, 
where $\mathcal{SC}$ and $\mathcal{BP}$ are symbols corresponding 
to \emph{separating circles} and \emph{bounding pairs} respectively. 
This function assigns to a vertex the symbol corresponding to its type. 
Namely, the symbol $\mathcal{SC}$ is assigned to isotopy classes of separating circles,
and the symbol $\mathcal{BP}$ to bounding pairs.

Alternativly, one can say simply that 
$\mathcal{SC}$-vertices are marked by $\mathcal{SC}$,
and $\mathcal{BP}$-vertices are marked by $\mathcal{BP}$.

\mytitle{Marked triangles.} Consider an unordered triple $C_0, C_1, C_2$ 
of disjoint non-separating cirlces in $S$.
If circles $C_0, C_1, C_2$ are pair-wise non-isotopic 
and pairs $\{C_0, C_1\}$, $\{C_1, C_2\}$, $\{C_2, C_0\}$ 
are bounding pairs of circles in $S$, then the
triple $C_0, C_1, C_2$ is called a \emph{marked triples of circles}. 
 
For a marked triple of circles $C_0, C_1, C_2$ 
the pairs $\{\vv{C_0},\vv{C_1}\}$, $\{\vv{C_1},\vv{C_2}\}$, $\{\vv{C_2},\vv{C_0}\}$
are bounding pairs, and hence $\mathcal{BP}$-vertices of Torelli building $\tts$.\sss 
Since the circles in the set $\mbox{V}{\eeq}\{C_0,C_1,C_2\}$ are pair-wise disjoint,
these three $\mathcal{BP}$-vertices are vertices of a simplex $\Delta$ of $\tts$.
In addition, these $\mathcal BP$-vertices are all different,
since circles $C_0, C_1, C_2$ are pair-wise non-isotopic. 
It follows that $\Delta$ is a simplex of dimension $2$, i.e. a triangle.

Any triangle of $\tts$ resulting from this construction is
called a \emph{marked triangle},\/ and there are no other marked triangle.
This designation of some triangles in $\tts$ as marked triangles is the
second part of our additional structure of $\tts$. 

\mytitle{A description of marked triples of circles.} Let $C_0, C_1, C_2$ be a marked triple of circles.
Consider two different elements $i,j$ of the set $\{0,1,2\}$. Let $k$ be the third element of this set.
Then $C_i\cup C_j$ divides $S$ into two parts, one of which contains $C_k$, and the other does not.
Let us denote by $S_{i,j}$ the part not containing $C_k$.
 
The union of subsurfaces $S_{0,1}$, $S_{1,2}$, and $S_{2,0}$ is a closed subsurface of $S$.
Since $S$ is a connected closed surface, this implies that $S$ is equal to this union, i.e.
\[
S {\eeq} S_{0,1} \cup S_{1,2} \cup S_{2,0}.
\]
The pair-wise intersections of our subsurfaces are $S_{0,1}\cap S_{1,2} {\eeq} C_1$, 
$S_{1,2}\cap S_{2,0} {\eeq} C_2$, and $S_{2,0}\cap S_{0,1} {\eeq} C_0$.
Each subsurface $S_{0,1}$, $S_{1,2}$, $S_{2,0}$ is a surface with exactly two boundary components.
Since the boundary components of each of these subsurfaces are non-isotopic in $S$ 
by the definition of marked triples of circles, none of these subsurfaces is an annulus.
This implies that each of them is a surface of genus $\geq 1$. 

This description implies, in particular, that $\gen{\eeq}\genus(S)\geqs 4$
if there is a marked triangle in $\tts$. 
Since most of our arguments depend on the existence of marked triangles, 
we usually have to assume that $\genus(S)\geqs 4$.

\mypar{The Torelli buildings.}{torelli-geom} The simplicial complex $\tts$ together with its marking of vertices
and the set of marked triangles is called the \emph{Torelli building\/} of $S$.
The Torelli building of $S$ is intended to play in the theory of Torelli groups a role 
similar to the role of Tits buildings in the theory of arithmetic groups 
and the role of the complexes of curves in the theory of Teichm{\"u}ller modular groups. 
Our main result about $\tts$ is a complete description of automorphisms of $\tts$; 
see Theorem \ref{automorphisms of geometry} below.

Recall (see the subsection \ref{cc}) that every diffeomorphism $F\colon S\rightarrow S$ induces 
a bijective self-map $\ffc$ of set of vertices of $\ccs$. 
Recall that every vertex of $\tts$ is either a $\mathcal{SC}$-vertex, or a $\mathcal{BP}$-vertex,
that $\mathcal{SC}$-vertices are at the same time (some) vertices of $\ccs$,
and that $\mathcal{BP}$-vertices are some pairs of vertices of $\ccs$.
Since $\ffc(\aavv{C}){\eeq}\aavv{F(C)}$ for every non-trivial circle $C$ in $S$,
the map $\ffc$ takes $\mathcal{SC}$-vertices of $\tts$ to $\mathcal{SC}$-vertices,
and takes pairs of vertices of $\ccs$ which are $\mathcal{BP}$-vertices of $\tts$ 
into $\mathcal{BP}$-vertices.
Therefore $\ffc$ induces a self-map of the set of the vertices of $\tts$. 
Obviously, this self-map is a bijection and, moreover, 
is an automorphism of the Torelli building $\tts$.
We will denote this automorphism by $\fft\colon\tts\tto\tts$.
It turns out that all automophisms of $\tts$ have the form $\fft$, at least if $\genus(S)\geq 5$.

\mypar{Main Theorem.}{automorphisms of geometry} \emph{If\/ $\gen{\eeq}\genus(S)\geqs 5$, then
every automorphism of the Torelli building\/ $\tts$ has the form\/ $\fft\colon\tts\tto\tts$,
where\/ $F\colon S{\tto} S$ is a diffeomorphism.}

Theorem \ref{automorphisms of geometry} 
is proved by showing that for $g\geqs 5$ every automorphism
of the Torelli building $\tts$ canonically induces an automorphism
of the complex of curves $C(S)$, and invoking
the theorem of the second author to the effect that all automorphisms
of $\ccs$ are induced by diffeomorphisms of $S$. 
At the same time Theorem \ref{automorphisms of geometry} is a Torelli 
analogue of this theorem about automorphisms of $\ccs$.

A key role in the proof of Theorem \ref{automorphisms of geometry} is played 
by the connectedness properties of 
two other simplicial complexes. 
The first one is the simplicial complex $T_{\rm sep}(R)$ 
defined only for connected surfaces $R$ with two boundary components. 
For such a surface $R$, the complex $T_{\rm sep}(R)$ has
as vertices the isotopy classes of circles separating two boundary components
of $R$ from each other, i.e. the isotopy classes of circles $D$ such that 
two components of $\partial R$ are contained in different components of $R\llsetmiss D$. 
The simplices of $T_{\rm sep}(R)$ are defined in the same way as the simplices of 
the complex of curves.

\mypar{Theorem \textup{(}Connectedness of $T_{\rm sep}(R)$\textup{)}.}{special connectedness} \emph{Let\/ $R$ 
be a connected surface with\/ $2$ boundary components. If\/
$\genus(R) \geqs 3$, then the complex\/ $T_{\rm sep}(R)$ is connected.}

The second complex defined for all connected surfaces $R$.
Let $R$ be such a surface, 
and let us define a \emph{genus $1$ circle}\/ in $R$
as a non-trivial separating circle in $R$ bounding in $R$ a torus with one hole in $R$. 
The \emph{complex of genus $1$ curves}\/ is defined as the full
subcomplex of  
$\ccc(R)$ spanning the set of the isotopy classes of genus
$1$ circles. It is denoted by $T_1(R)$.

\mypar{Theorem \textup{(}Connectedness of\/  $T_1(R)$\textup{)}.}{connectedness} \emph{Let\/ $R$ 
be a compact connected surface. If\/ $\genus(R)\geqs 3$, 
then the  complex $T_1(R)$  is connected.} 

Both connectedness theorems do not hold without the assumption $\genus(R)\geqs 3$;
in fact, if $\genus(R){\eeq}2$, then complexes $T_{\rm sep}(R)$
and $T_1(R)$ 
have an infinite number of vertices and no edges,
and if $\genus(R)\leq 1$, these complexes are empty, i.e. have no vertices either.

\mysection{Proof of connectedness theorems}{proof-conn-ths}

We will prove Theorem \ref{connectedness} first (see the subsection \ref{proof-connecteness}), and 
then deduce Theorem \ref{special connectedness} from it (see the subsection \ref{proof-special-connectness}). 
This deduction will be based on ideas similar to the ideas of the proof of Theorem~\ref{connectedness}.

\mypar{Auxilary complexes and their connectedness.}{aux}
Our proofs of Theorems \ref{connectedness}, \ref{special connectedness} 
are based on connectedness of some auxilary complexes, discussed in this subsection.
Let $R$ be a compact connected surface, possibly with non-empty
boundary, in contrast with $S$.

\mytitle{Complex $ \ccc_0(R)$.} The first one is the subcomplex $ \ccc_0(R)$ 
of the complex of curves $ \ccc(R)$ having as vertices the isotopy classes of non-separating
circles and as simplices collections of isotopy classes of several
disjoint circles which are jointly nonseparating, 
i.~e. such that the complement of their union in $R$ is connected.

\mylemma{Lemma.}{c-0} \emph{The complex $ \ccc_0(R)$ is connected if the genus of $R$ is $\geqs 3$.}

This result was proved independently by J. Harer \cite{h1, h2} and
the second author \cite{ivcd}. J. Harer and the second author also proved
and used various versions of this result; see \cite{h1, h2}, \cite{i-usp,
i-stab, Iv1}. The proof proceeds by proving the connectedness of the
complex of curves $C(S)$ first and then using some combinatorial arguments.

\mytitle{Complexes $\mathcal{H}(R; y_1, y_2)$ and $\mathcal{H}_0(R; y_1, y_2)$.} Suppose that $R$ is a 
connected surface with at least $2$ boundary components and that $y_1$, $y_2$
are $2$ points belonging to \emph{different}\/ components of 
the boundary $\partial R$. Both complexes $\mathcal{H}(R; y_1, y_2)$ and 
$\mathcal{H}_0(R; y_1, y_2)$ have as vertices the isotopy classes of embedded arcs connecting
$y_1$ with $y_2$; the isotopies are assumed to be point-wise fixed on the
boundary $\partial R$. The simplices of $\mathcal{H}(R; y_1, y_2)$ 
are collections of isotopy classes of arcs with disjoint interiors. 
The simplices of $\mathcal{H}_0(R; y_1, y_2)$ are collections of isotopy
classes of arcs with disjoint interiors and \emph{jointly non-separating}\/ in $R$.
These complexes were introduced by J. Harer \cite{h1} and was used by him and 
by the second author (see, for example, \cite{i-stab,Iv1})
in the proofs of homological stability theorems for groups $\Mod_R$.

The complex $\mathcal{H}(R; y_1, y_2)$ is introduced here in order to place
the following lemma in its proper context and to discuss 
the somewhat subtle situation with its proof (see remarks
after the lemma); it will not be used in this paper directly.

\mylemma{Lemma.}{h-h-conn} \emph{The complex\/ $\mathcal{H}_0(R; y_1, y_2)$
is connected if the genus of\/ $R$ is\/ $\geqs 2$.}

This is a special case of Theorem 1.4 in \cite{h1}. The proof proceeds by
proving first that the complexes $\mathcal{H}(R; y_1, y_2)$ are highly connected; 
this is the content of Theorem 1.6 in \cite{h1}. 
The connectivity proprerties of $\mathcal{H}(R; y_1, y_2)$ imply (weaker)
connectivity properties of $\mathcal{H}_0(R; y_1, y_2)$
by a standard by now argument of Harer 
(see \cite{h1}, proofs of Theorems 1.1, 1.3 and 1.4). 

The proof of Theorem 1.6 in \cite{h1} turned out to be flawed, 
as was discovered by N. Wahl \cite{wahl-non-or} in connection with an investigation of 
homological stability for mapping class groups of non-orientable surfaces. 
Namely, Harer \cite{h1} does not treat the case of two points belonging to the same component of boundary.
While we do not need this case of Harer's theorem, Harer's proof proceeds by an inductive argument in which
both cases are needed to carry on the induction. 

A correct proof of Theorem 1.6 from \cite{h1}  
in the full generality was provided by N. Wahl \cite{wahl-non-or}; 
see \cite{wahl-non-or}, Theorem 2.3. 

\mytitle{Complexes $H(R,C,x)$ and $H_0(R,C,x)$.} These complexes are, in fact, versions
of complexes $\mathcal{H}(R; y_1, y_2)$ and $\mathcal{H}_0(R; y_1, y_2)$. 
Suppose that $R$ is a connected surface and $C$ is a non-separating circle on $R$. 
Suppose that a point $x\in C$ is fixed. Both complexes have as vertices
the isotopy classes of circles $D$ on $R$ such that $C\cup D {\eeq}\{x\}$
and $D$ intersects $C$ (at the point $x$) transversely; 
the isotopies are assumed to be pointwise fixed on $C$.
The simplices of $\mathcal{H}(R,C,x)$ 
are collections of isotopy classes of circles $D_1$, $D_2$, \ldots, $D_n$ such that
the complements $D_1\setmiss x$, $D_2\setmiss x$, \ldots, $D_n\setmiss x$ are disjoint. 
The simplices of $\mathcal{H}_0(R,C,x)$ are collections of isotopy classes of arcs such that
the complements $D_1\setmiss x$, $D_2\setmiss x$, \ldots, $D_n\setmiss x$ are disjoint, 
and the complement of $C\cup D_1\cup D_2\cup\ldots\cup D_n$ in $R$ is connected.

\mylemma{Lemma.}{h-c-conn} \emph{The complex\/ $\mathcal{H}_0(R,C,x)$
is connected if the genus of\/ $R$ is\, $\geqs 3$.}

\emph{Proof.} Let us cut our surface $R$ along the circle $C$. 
Let $R_C$ be the resulting surface, and let $x_1$, $x_2$
be two points on the boundary $\partial R_C$ resulting from the point
$x$. Clearly, the points $x_1$, $x_2$ belong to two different components
of $\partial R_C$ consists of two components. Moreover, 
the complex $H_0(R,C,x)$ is obviously isomorphic to $H_0(R_C; x_1, x_2)$ 
and the genus of $R_C$ is at least $2$ if the genus of $R$ is at least $3$.
Hence the lemma follows from Lemma \ref{h-h-conn}. \eproof

\mypar{Proof of Theorem \ref{connectedness}.}{proof-connecteness} Recall that $R$ is a compact connected orientable surface. 
As is well known and obvious, for every genus $1$ circle $C$ there are two circles $D$, $E$ 
transversely intersecting at a single point such that $C$ is isotopic to the boundary
of a regular neighborhood of the union $D\cup E$, or, what is the same,
such that $D\cup E$ is contained in the torus with one hole bounded by
$C$.  The circle $C$ is determined up to isotopy by circles $D$ and $E$;
slightly abusing the notations, we will write $C{\eeq}C(D,E)$. Note that for a given circle $C$
such a pair $D$, $E$ is far from unique even up to isotopy.

Let $C$, $C'$ be two genus $1$ circles which we would like to
connect by a sequence of edges of the complex of genus $1$ circles. 
Let us choose pairs of circles $D$, $E$ and $D'$, $E'$ such that $C{\eeq}C(D,E)$
and $C'{\eeq}C(D',E')$. Notice that circles $D$, $D'$ are non-separating.
By Lemma \ref{c-0} the complex $ \ccc_0(R)$ is connected, and hence there is a sequence of non-separating
circles 
\[
\seq{D}{D}{D'}{n}
\] 
such that $D_i$ is connected with $D_{i+1}$ by an edge of $ \ccc_0(R)$ for each $i$.
In other words, for each $i$, the circles $D_i$, $D_{i+1}$ are isotopic to two jointly non-separating disjoint circles.  
In fact, after replacing the circles $D_i$ by isotopic circles, if necessary (for example, one can
fix a hyperbolic metric with geodesic boundary on $R$ and take geodesic circles isotopic to
the circles $D_i$), we may assume that the circles $D_i$, $D_{i+1}$ are themselves disjoint
and jointly non-separating, i.e. the complement of $D_i\cup D_{i+1}$ is connected. 

Since the complements  $R\setminus (D_i\cup D_{i+1})$ are connected, we can choose for each $i{\eeq}1,\ldots$, $n\minus 1$
a circle $E_i$ transversely intersecting $D_i$ at a single point and disjoint from $D_{i+1}$. Let $E_n{\eeq}E'$,
and let us assume that $E_1{\eeq}E$. Note that the circles $E_i$ and $E_{i+1}$ may be not disjoint, even up to
isotopy. We will deal with this difficulty by introducing some additional curcles $F_{i+1}$.

Note that the union $D_i\cup E_i$ is disjoint from $D_{i+1}$, and hence a regular neighborhood $N_i$
of $D_i\cup E_i$ is disjoint from $D_{i+1}$. Since the complement of $D_i\cup D_{i+1}$ is connected,
the result $R_{\rm cut}$ of cutting $R$ along the union $D_i\cup D_{i+1}$ is also connected. After the cutting,
$E_i$ turns into an arc $A_i$ connecting two different components of the boundary $\partial R_{\rm cut}$. 
Since $R_{\rm cut}$ is connected, the complement of this arc in $R_{\rm cut}$ is also connected. 
After the cutting, the regular neighborhood $N_i$ turns into a 
regular neighborhood of the union of $A_i$ with these two boundary components. 
It follows that the complement of $N_i$ in $R_{\rm cut}$ is connected.
Moreover, the genus of $R_{\rm cut}\setmiss N_i$ is equal to the genus of $R_{\rm cut}$.
Since the genus of $R$ is\, $\geqs 3$, the genus of $R_{\rm cut}$ is\, $\geqs 1$, 
and hence the genus of $R\llsetmiss(N_i\cup D_{i+1})$ is also\, $\geqs 1$. 

The last statement implies that there is a circle $F_{i+1}$ transversely intersecting $D_{i+1}$ 
at a single point and disjoint from $N_i$. In particular,
$F_{i+1}$ is disjoint from $E_i$. We may assume that, in addition, 
$F_{i+1}$ intersects $D_{i+1}$ at the same point as $E_{i+1}$. 

Clearly, $C{\eeq}C(D,E){\eeq}C(D_1 ,E_1)$, $C'{\eeq}C(D',E'){\eeq}C(D_n ,E_n)$, 
and for every $i{\eeq}1,\ldots$, $n{\minus}1$ the circles $C(D_i, E_i)$ and $C(D_{i+1}, F_{i+1})$ are disjoint up to isotopy.  
If for every $i{\eeq}1,\ldots$, $n{\minus}1$ we will be able to connect the isotopy classes 
$\langles C(D_{i+1}, F_{i+1})\rangles $ and $\langles C(D_{i+1}, E_{i+1})\rangles$ by a chain of edges in the complex of genus $1$ circles, 
we will be able to solve our problem: the chain formed by union of such chains and edges connecting
$\langles C(D_i, E_i)\rangles$ with $\langles C(D_{i+1}, F_{i+1})\rangles$ will connect the isotopy classes $\langles C\rangles$ and $\langles C'\rangles$.

So, let us fix $i$, $1\leqs i \leqs n{\minus}1$. Curves $F_{i+1}$ and $E_{i+1}$
define two vertices of the complex  $H_0(R,D_{i+1},x)$, where $x$ is
the point of intersection of $E_{i+1}$ and $F_{i+1}$ with $D_{i+1}$.
By Lemma \ref{h-c-conn} the complex $H_0(R,D_{i+1},x)$ is connected.
Hence we can find a sequence of
circles  
\[
\seq{G}{F_{i+1}}{E_{i+1}}{m}
\] 
such that the circles $G_j$, $G_{j+1}$ are isotopic by an isotopy fixed on $D_{i+1}$ 
to two disjoint circles $G'_j$, $G'_{j+1}$ such that 
the complement of union $D_{i+1}\cup G'_{j} \cup G'_{j+1}$ in $R$ does not separate $R$
for $1\leqs j \leqs m-1$. As usual, we may assume that curves $G_j$, $G_{j+1}$
are themselves disjoint and the complement of union 
$D_{i+1}\cup G_{j} \cup G_{j+1}$ in $R$ does not separate $R$. 

Consider now the regular neighborhood of the union $D_{i+1}\cup G_{j} \cup
G_{j+1}$. Since the regular neighborhood of $D_{i+1}\cup G_{j}$ is a
torus with one hole, and the regular neighborhood of $D_{i+1}\cup G_{j}
\cup G_{j+1}$ results from it by attaching a band to the boundary, the
latter regular neighborhood is a torus with two holes.  This torus with
two holes has a connected complement in $R$ because $D_{i+1}\cup G_{j}
\cup G_{j+1}$ does not separate $R$. Since the genus of $R$ is 
$\geqs 3$, this complement contains a genus $1$ circle $H_j$. Now we can
connect $\langles C(D_{i+1}, G_j)\rangles$ and $\langles C(D_{i+1},
G_{j+1})\rangles$ by a chain consisting of two edges: one connecting
$\langles C(D_{i+1}, G_j)\rangles$ with $\langles H_j\rangles$ and the other
connecting $\langles H_j\rangles$ with $\langles C(D_{i+1},
G_{j+1})\rangles$. By taking all these two-edge chains together we will
get a chain connecting $\langles C(D_{i+1}, F_{i+1})\rangles$ with
$\langles C(D_{i+1}, E_{i+1})\rangles$.  As was explained in the
previous paragraph, this implies that we can connect $\langles C\rangles$
with $\langles C'\rangles$. This completes the proof of 
Theorem \ref{connectedness}. \eproof

\mypar{Proof of Theorem \ref{special connectedness}.}{proof-special-connectness}
Recall that in this theorem $R$ is a compact connected surface with exactly $2$ boundary components.
As before, we will call a circle $C$ on $R$ a \emph{genus $1$ circle}\/ 
if $C$ bounds in $R$ a torus with one hole. 
Let $B_1$ and $B_2$ be the components of $\partial R$.  
Let us call a circle $C$ in $R$ \emph{special}\/ if $C$ bounds together
with either $B_1$ or $B_2$ a torus with two holes. 
Since every circle separating $B_1$ from 
$B_2$ is disjoint from some special circle, it sufficient to prove
the connectedness of full subcomplex $T_{\rm spec}(R)$ of $T_{\rm sep}(R)$ having the set of the 
isotopy classes special circles as the set of vertices.
We will call the complex $T_{\rm spec}(R)$ the \emph{complex of special circles}\/ in $R$.

For a genus $1$ circle $D$ on $R$, we will denote by $\tau_D$ be the torus with one
hole bounded by $D$ in $R$. Suppose that $J$ is an embedded arc in $R$ connecting $\tau_D$
with boundary $\partial R$ and intersecting $\tau_D$ and $\partial R$ only at its endpoints.
Let $B_J{\eeq}B_1\mbox{ or }B_2$ be the component of $\partial R$ connected with $\tau_D$ by $J$.   
Then the boundary of a regular neighborhood
of the union $\tau_D\cup J\cup B_J$ is clearly a special circle. By a
slight abuse of notation, we will denote this circle by $C(D,J)$. It
is clear that every special circle has the form $C(D,J)$ for some $D$, $J$.
The remaining part of the proof is parallel to the proof of Theorem
\ref{connectedness}, with the arcs $J$ playing a role similar to the
circles $E_i$ in that proof.

Suppose that $C$, $C'$ are two special circles which we would like to
connect by a chain of edges in $T_{\rm spec}(R)$. Since
every special circle bounding together with $B_2$ a torus with two holes is disjoint
from some special curve bounding a torus with two holes together with $B_1$, we may
assume that both curves $C$ and $C'$ bound a torus with two holes together with $B_1$ .  

Let us begin by choosing some genus $1$ curves $D$, $D'$ and arcs $J$, $J'$ 
such that $C{\eeq}C(D,J)$ and $C{\eeq}C(D',J')$. By Theorem \ref{connectedness} the complex $T_1(R)$ of genus $1$
circles is connected. Hence there is a sequence of genus $1$ circles 
\begin{equation}
\label{chain}
\seq{D}{D}{D'}{n}
\end{equation}
such that $D_i$ is connected with $D_{i+1}$ by an edge of $T_1(R)$ for every $i{\eeq}1,2,\ldots, n{\minus}1$.
In other words, for every $i{\eeq}1,2,\ldots, n{\minus}1$ the circles $D_i$, $D_{i+1}$ are isotopic to two disjoint circles. 
As usual, we can assume that circles $D_i$, $D_{i+1}$ are themselves disjoint. 
Let us choose an arc $J_i$ connecting $\tau_{D_i}$ with $B_1$ for every $i{\eeq}1,2,\ldots, n$
(in particular, $B_{J_i}{\eeq}B_1$ for every $i{\eeq}1,2,\ldots, n$). 
We may assume that $J_1{\eeq}J$ and $J_n{\eeq}J'$. 

Let us fix for a moment some $i{\eeq}1,2,\ldots,n{\minus}1$.
Since the circles $D_i$ and $D_{i+1}$ are disjoint, either
subsurfaces $\tau_{D_i}$ and $\tau_{D_{i+1}}$ are disjoint, or one of them 
is contained in the other (because each of them is a torus with one hole). 
In the second case their boundaries $D_i$\/, $D_{i+1}$ are isotopic and 
the corresponding vertices of $T_1(R)$ are equal.
In this case we can shorten our chain (\ref{chain}) by removing either $D_i$ or $D_{i+1}$.
Therefore, we can assume that subsurfaces $\tau_{D_i}$ and $\tau_{D_{i+1}}$ are disjoint
for every $i{\eeq}1,2,\ldots,n{\minus}1$.

The rest of the proof would be easier if we would be able 
to assume that for every $i{\eeq}1,2,\ldots,n{\minus}1$ 
the arcs $J_i$ and $J_{i+1}$ are disjoint at least up to an isotopy 
(of course, only the isotopies fixed on $B_1$ and on $\tau_{D_i}$\/, $\tau_{D_{i+1}}$ 
respectively are allowed). But this is not possible in general. 
In order to deal with this difficulty, 
let us choose for every $i{\eeq}1,2,\ldots, n{\minus}1$ an arc $A_{i+1}$ connecting 
$\tau_{D_{i+1}}$ with $B_1$ and disjoint from $\tau_{D_i}\cup J_i\cup B_1$
(such an arc exists because now we are assuming that 
$\tau_{D_i}$ and $\tau_{D_{i+1}}$ are disjoint).

\textsc{Claim.} \emph{For every $i{\eeq}1,2,\ldots,n{\minus}1$ the vertices
$\langles C(D_i\/,J_i) \rangles$ and $\langles C(D_{i+1}\/,A_{i+1}) \rangles$ 
can be connected in $T_{\rm spec}(R)$ by a chain consisting of no more than $2$\/ edges.}

\emph{Proof of the claim.} As we noted above, every regular neighborhood of 
$\tau_{D_i}\cup J_i \cup B_1$ is a torus with one hole.  
We can choose a regular neighborhood $N_i$ of $\tau_{D_i}\cup J_i \cup B_1$ such
that it is disjoint from $\tau_{D_{i+1}}$ and that its intersection with the arc $A_{i+1}$ 
is a segment contained in $A_{i+1}$ and having one endpoint in $B_1$ and the other endpoint in $C(D_i\/, J_i)$. 
The we can construct a regular neighborhood $L_i$ of 
\[
\tau_{D_i}\cup J_i \cup B_1 \cup \tau_{D_{i+1}} \cup A_{i+1}
\] 
by adding to $N_i$ a regular neighborhood $\nu_{D_{i+1}}$ of $\tau_{D_{i+1}}$ 
(which is a subsurface of $R$ containing the subsurface $\tau_{D_{i+1}}$ and isotopic to it)
and a band connecting $N_i$ with $\nu_{D_{i+1}}$.
Since both surfaces $N_i$ and $\nu_{D_{i+1}}$ are surfaces of genus $1$,
the genus of $L_i$ is $2$. 

Since the genus of $R$ is assumed to be $\geqs 3$,
there is a special curve $H_i$ bounding a torus with two holes
together with the component $B_2$ of $\partial R$ and
disjoint from $L_i$. In particular, $H_i$ is disjoint from
the union $\tau_{D_i}\cup J_i \cup \tau_{D_{i+1}}\cup A_{i+1}\subset L_i$.
It follows that $\langles C(D_i,J_i) \rangles$ can be connected with
$\langles C(D_{i+1}\/,A_{i+1}) \rangles$ by a chain consisting of two edges:
one connecting $\langles C(D_{i}\/,J_{i})\rangles$ with $\langles H_i\rangles$
and the other connecting $\langles H_i\rangle$ with $\langles C(D_{i+1}\/,A_{i+1})\rangles$.
This proves the claim. \esubproof

It remains to connect $\langles C(D_{i+1}\/,A_{i+1}) \rangles$ with 
$\langles C(D_{i+1}\/,J_{i+1}) \rangles$ by a chain in $T_{\rm spec}(R)$ for all $1\leqs i\leqs n{\minus}1$. 
By using, if necessary, an isotopy of $R$ preserving $\tau_{D_{i+1}}$ set-wise 
(and not necessarily fixed on the boundary), we can assume that $J_{i+1}$
and $A_{i+1}$ have the same endpoints, which we denote by $x_1$,
$x_2$.  Let $R_{i+1}{\eeq}(R\llsetmiss\tau_{D_{i+1}})\cup D_{i+1}$. 
Consider the complex $H_0(R_{i+1}; x_1,x_2)$.  
Since the genus of $R$ is assumed to be $\geqs 3$, the genus of $R_{i+1}$ is $\geqs 2$, 
and hence the complex $H_0(R_{i+1}; x_1,x_2)$ is connected by Lemma \ref{h-h-conn}.  
Therefore there is a sequence of arcs
\[
\seq{K}{J_{i+1}}{A_{i+1}}{m}
\]
such that $K_j$ is connected with $K_{j+1}$ by an edge of $H_0(R_{i+1}; x_1,x_2)$ for every $j{\eeq}1,2\ldots,m-1$.
In other words, for every $j{\eeq}1,2,\ldots, m-1$ the arcs $K_j$, $K_{j+1}$ 
are isotopic in $R_{i+1}$ to two disjoint arcs
by isotopies point-wise fixed on $\partial R_{i+1}$.
As usual, we can assume that the interiors of $K_j$, $K_{j+1}$ are
disjoint for every $j{\eeq}1,2,\ldots, m-1$. Then, in addition, $K_j\cup K_{j+1}$
does not separates $R_{i+1}$ because $\langles K_j\rangles$ is connected
with $\langles K_{j+1}\rangles$ by an edge of $H_0(R_{i+1}; x_1,x_2)$. 

Since the genus $R_{i+1}$ is $\geqs 2$ and $K_j\cup K_{j+1}$ does not separates $R_{i+1}$, 
for every $j{\eeq}1,2,\ldots$, $m-1$ the surface $R_{i+1}\setmiss(K_j\cup K_{j+1})$ has genus $\geqs 1$. 
Therefore $R_{i+1}\setmiss (K_j\cup K_{j+1})$ contains a circle $G_j$ 
bounding together with $B_2$ a torus with two holes.
Clearly, $G_j$ is disjoint from the union $K_j\cup K_{j+1}$.
Considered as a circle on $R$, this circle $G_j$ is disjoint from the
union $\tau_{D_{i+1}}\cup K_j\cup K_{j+1}$. It follows that $\langles
C(D_{i+1},K_j) \rangles$ can be connected with $\langles C(D_{i+1},K_{j+1}) \rangles$ 
in $T_{\rm spec}(R)$ by a chain consisting
of two edges: the edge connecting $\langles C(D_{i+1},K_j) \rangles$ with
$\langles G_j\rangles$ and the edge connecting $\langles G_j\rangles$ with
$\langles C(D_{i+1},K_{j+1}) \rangles$.  By joining all these two-edge
chains together, we get a chain connecting $\langles C(D_{i+1},A_{i+1})
\rangles$ with $\langles C(D_{i+1},J_{i+1}) \rangles$.  As we noticed
above, this completes the proof of Theorem \ref{special connectedness}. \eproof
\vspace{4ex}

\mysection{Encodings}{encodings}

This section is devoted to a 
key tool need for the proof Theorem \ref{automorphisms of geometry},
namely, to encodings of isotopy classes of 
non-separating circles in $S$ in terms of the Torelli building $\tts$.

\mypar{Encodings of non-separating circles.}{en-nosep} 
Let $\gamma$ be a vertex of $\tts$ represented by a bounding pair of
non-separating circles $C_0$, $C_1$. We would like to find a description of
the isotopy class $\langles C_0 \rangles$ of the circle $C_0$ alone
in terms of the Torelli building $\tts$. 
The idea is to use a separating circle $D$ on $S$ 
which is up to isotopy disjoint from $C_1$
and which cannot be made by an isotopy disjoint from $C_0$.
In terms of the complex of curves $C(S)$, 
the isotopy class $\gamma_0{\eeq}\langles C_0 \rangles$ is a vertex of $C(S)$, 
and our condition means that  the isotopy class $\delta{\eeq}\langles D \rangles$ 
(which is also a vertex of $C(S)$) 
is connected by an edge of $C(S)$ with $\langles C_1 \rangles$ 
and is not connected by an edge of $C(S)$ with $\langles C_0 \rangles$.
  
Obviously, the vertex $\langle C_0\rangle$ of $C(S)$ can be recovered from 
the \emph{ordered}\/\/ pair $(\gamma, \delta)$, where $\delta{\eeq}\langle D\rangle$.
The ordered pair $(\gamma, \delta)$ can be considered as an \emph{encoding}\/ of 
the vertex $\gamma_0$. Clearly, there are many encodings of $\gamma_0$,
and we need to be able to tell when two pairs encode the same vertex. 
In addition, we need to be able to tell when a pair is an encoding
of some non-separating circle.

\mypar{Encoding by admissible pairs.}{encoding-add}
For technical reasons we will use only some pairs $(\gamma, \delta)$ as encodings,
namely, only the pairs which we will call \emph{admissible}\/ pairs. They are defined as follows.

\mytitle{Admissible pairs.} 
Let $\gamma$ be a vertex of $\tts$ corresponding to a bounding pair of circles and 
let $\delta$ be a vertex of $\tts$ corresponding to a separating circle.  
We call the pair $(\gamma, \delta)$ an \emph{admissible pair}\/ if 
there are two vertices $\gamma'$, $\beta$ of $\tts$ corresponding to bounding pair of circles
(i.e. marked by the symbol $\mathcal{BP}$ in the Torelli building $\tts$)
such that the following two conditions hold.
\vspace{-3ex}
\begin{enumerate}
\item[(i)] $\{ \gamma , \gamma' ,\beta\}$ is a marked triangle of $\tts$.
\vspace{-1ex} 
\item[(ii)] $\delta$ is connected by an edge of $\tts$ with $\beta$.
\vspace{-1ex}
\item[(iii)] $\delta$ is \emph{not}\/ connected by an edge with either $\gamma$ or $\gamma'$.
\end{enumerate}
\vspace{-3ex}
If the the conditions (i), (ii), and (iii) hold, we say that the (ordered) triple $(\gamma , \gamma' ,\beta)$ 
\emph{justifies the admissibility}\/ of the pair $(\gamma, \delta)$. 
Note that these conditions are stated entirely in terms of the Torreli geometry structure of $\tts$.

Since the definition of an admissible pair requires the existence of a marked triangle,
admissible pairs exist only if $\genus(S)\geq 4$ (see the end of the subsection \ref{add-structure}).

\mylemma{Lemma.}{singling out} \emph{Let $C_0$, $C_1$, $D$, and $\gamma$, $\delta$\/ be three circles on\/ $S$ and
two vertices of\/ $\tts$ with the properties described in the section\/ \emph{\ref{en-nosep}}. 
Suppose that the pair\/ $(\gamma,\delta)$ is admissible, and let $(\gamma, \gamma', \beta)$
be any triple justifying the admissibility of\/ $(\gamma,\delta)$. 
Then the isotopy class\/ $\gamma_0{\eeq}\langles C_0 \rangles$ of\/ $C_0$ is 
the unique element of the set\/ $\gamma\cap\gamma'{\eeq}\gamma\llsetmiss\beta{\eeq}\gamma'\llsetmiss\beta$, 
and $(\gamma,\delta)$ is an encoding of\/ $\gamma_0$.}

\proof Let $(\gamma, \gamma',\beta)$ be some triple justifying the admissibility of $(\gamma,\delta)$.
Then $\{ \gamma , \gamma' ,\beta\}$ is a marked triangle and is determined by
three circles as in the definition of the marked triangles in section \ref{tg-def}.
Since $\gamma$ is the vertex corresponding to the bounding pair of circles $C_0$, $C_1$,
we can assume that two of these three circles are $C_0$, $C_1$. Let $C_2$ be the third circle.

Then $\beta$ is the vertex of $\tts$ corresponding to the bounding pair of circles $C_i$, $C_j$,
where $i$, $j$ is a pair of two distinct elements of the set $\{0,1,2\}$. 
Since $\delta$ is not connected by an edge of $C(S)$ with $\gamma_0{\eeq}\langles C_0 \rangles$, 
the circle $D$ is not isotopic to a circle disjoint from $C_0$.
On the other hand $\delta$ is connected by an edge of $\tts$ with $\beta$ and
hence $D$ is isotopic to a circle disjoint from both $C_i$ and $C_j$.
By comparing the last two statements we conclude that neither $i$, nor $j$ could be
equal to $0$. It follows that $\{i,j\}{\eeq}\{1,2\}$, $\{C_i, C_j\}{\eeq}\{C_1, C_2\}$,
and $\beta$ is the vertex of $\tts$ corresponding to the bounding pair of circles $C_1$, $C_2$.

In other words, $\beta$ is equal to the bounding pair $\{\langles C_1 \rangles, \langles C_2 \rangles\}$.
Since $\gamma$ corresponds to the bounding pair of circles $C_0$, $C_1$, 
it is equal to the bounding pair $\{\langles C_0 \rangles, \langles C_1 \rangles\}$.
Since $\{\gamma,\gamma',\beta\}$ is a marked triangle, 
$\gamma'$ is equal to the bounding pair $\{\langles C_0 \rangles, \langles C_1 \rangles\}$.
Taken together, these descriptions of $\gamma$, $\gamma'$, $\beta$ imply that
$\gamma\cap\gamma'{\eeq}\gamma\llsetmiss\beta{\eeq}\gamma'\llsetmiss\beta$ and
the isotopy class $\langles C_0 \rangles$ is the unique element of this set. 

It remains to prove that $(\gamma,\delta)$ is an encoding of $\gamma_0$.
Indeed, since $\{ \gamma , \gamma' ,\beta\}$ is a marked triangle and $\delta$ is not connected
by an edge of $\tts$ with $\beta$, the circle $D$ is isotopic to a circle disjoint from $C_1\cup C_2$.
In particular, $D$ is isotopic to a circle disjoint from $C_1$. 
On the other hand, $\delta$ is not connected by an edge with $\gamma$.
It follows that $D$ is not isotopic to a circle disjoint from $C_0$.
Hence $(\gamma,\delta)$ is indeed an encoding of $\gamma_0$.
This completes the proof of the lemma. \eproof

\mylemma{Corollary.}{pairs-are-enc} \emph{If\/ $(\gamma,\delta)$ is an admissible pair,
then\/ $(\gamma,\delta)$ is an encoding of the isotopy class of a non-separating circle.} \eproof

\mypar{Moves of admissible pairs.}{moves} We will prove that two admissible pairs
encode the same isotopy class of non-separating circles if and only if
these two pairs can be connected by a sequence of moves of two types. We start
with the definition of these moves.

\mytitle{Moves of type I.}\, If a triple $(\gamma, \gamma', \beta)$ justifies the admissibility
of the pair $(\gamma, \delta)$, then there is a \emph{move of type I}\/\/ from $(\gamma, \delta)$
to $(\gamma', \delta)$ (and there are no other moves of type I). 
Alternatively, we say such pairs $(\gamma, \delta)$ and $(\gamma', \delta)$
are \emph{connected by a move of type I}.

\mytitle{Moves of type II.}\,  If a triple $(\gamma, \gamma', \beta)$ simultaneously justifies the admissibility
of the pairs $(\gamma, \delta)$ and $(\gamma, \delta')$, then there is a \emph{move of type II}\/\/
replacing $(\gamma, \delta)$ by $(\gamma, \delta')$ (and there are no other moves of type II). 
Alternatively, we say such pairs $(\gamma, \delta)$ and $(\gamma, \delta')$
are \emph{connected by a move of type II}.

\mylemma{Lemma.}{not changed} \emph{The isotopy class of non-separating circles encoded by
an admissible pair is not changed by moves of both types.}

\emph{Proof.} We may assume that we are in the situation of
Lemma \ref{singling out}. Then by this lemma $\gamma_0$
is the unique element of $\gamma\cap\gamma'$. 
This immediately implies our claim for moves of type I. 
In addition, this implies that
if a triple justifying the admissibility of $(\gamma, \delta)$ 
is known, we can find $\gamma_0$ by using only this triple
without a recourse to $\delta$. 
This implies our claim for moves of type II. \eproof

\mypar{Theorem (Changing $\delta$).}{changing delta} \emph{Suppose that the genus of $S$ is $\geqs 5$.
Let $(\gamma, \delta')$ and $(\gamma, \delta'')$ be two admissible pairs with the
same first vertex and encoding the same isotopy class of non-separating
curve. Then there is a sequence of moves of the type II starting at
$(\gamma, \delta')$ and ending at $(\gamma, \delta'')$.}

\emph{Proof.} Let $(\gamma, \gamma',\beta')$ and 
$(\gamma, \gamma'',\beta'')$ be some triples justifying the
admissibility of $(\gamma, \delta')$ and $(\gamma, \delta'')$
respectively. Suppose that $C_0$, $C_1$ is a bounding pair of circles
such that $\gamma$ is equal to the corresponding bounding pair.
We can assume that $(\gamma, \delta')$ and $(\gamma, \delta'')$
encode the isotopy class $\gamma_0{\eeq}\langles C_0\rangles$ of $C_0$.
Since the triangles $\{ \gamma, \gamma',\beta'\}$ and 
$\{\gamma, \gamma'' ,\beta''\}$ are marked, they correspond to triples of 
circles $C_0$, $C_1$, $C'_2$ and $C_0$, $C_1$, $C''_2$ respectively 
for some non-separating circles $C'_2$ and $C''_2$. 

By Lemma \ref{singling out}, $\gamma_0$ is the unique element of
each of the intersections $\gamma\cap\gamma'$ and $\gamma\cap\gamma''$.
In view of the definition of marked triangles, this implies that
$\gamma'$, $\gamma''$ are the vertices of $\tts$ corresponding 
to the bounding pairs of circles $\{C_0,C'_2\}$ , $\{C_0,C''_2\}$ respectively.
Similarly, $\beta'$, $\beta''$ are the vertices of $\tts$ corresponding 
to the bounding pairs of circles $\{C_1,C'_2\}$ , $\{C_1,C''_2\}$ respectively.

Since $C_0$, $C_1$ is a bounding pair of circles, 
the set difference $S\setmiss(C_0\cup C_1)$ consists two components.
We will denote the closures of these components by $S'$ and $S''$.
Clearly, $S'$ and $S''$ are subsurfaces of $S$ with the common boundary $C_0\cup C_1$ 
The rest of the proof splits into two cases. 

\mytitle{Case 1. The circles\/ $C'_2$ and\/ $C''_2$ are contained in different components of\/ $S\setmiss(C_0\cup C_1)$.}
We may assume that $C'_2$ is contained $S'$ and $C''_2$ is contained in $S''$.
 
In this case $C_0$ and $C'_2$ bound a subsurface of $S'$ of genus at least $1$
with boundary equal to $C_0\cup C'_2$. 
Clearly, this subsurface is disjoint from $C_1$.
Its genus is $\geqs 1$, because otherwise it is an annulus
and $C_0$ and $C_2$ are isotopic, which is not allowed by the definition of marked triangles.
Similarly, $C_0$ and $C''_2$ bound a subsurface of $S''$ disjoint from $C_1$. 
Its genus is also $\geqs 1$ and its boundary is equal to $C_0\cup C''_2$.
These two subsurfaces are situated on different sides of $C_0$ because $C'_2$ and $C''_2$ are
situated on different sides of $C_0\cup C_1$ (i.e. are contained in different components of
$S\setmiss(C_0\cup C_1)$).
The union $Q$ of these two subsurfaces is a subsurface of $S$ disjoint from $C_1$. 
Its boundary is equal to $C'_2\cup C''_2$ and its genus is $\geqs 2$. 
Clearly, $C_0$ is a nontrivial separating circle in $Q$. Moreover, $C_0$
separates the components $C'_2$ and $C''_2$ of $\partial Q$.
As is well known, this implies that there is a nontrivial separating circle 
$D_1$ in $Q$ which is not isotopic in $Q$ to a circle disjoint from $C_0$.
One can obtain such a circle $D_1$ by applying to $C_0$ a Dehn twist
along any circle not isotopic to a circle disjoint from $C_0$. 
The point of this argument is that it is easier to imagine or draw a 
non-separating circle with this property than a separating one.

\textsc{Claim.} \emph{Let $\delta_1{\eeq}\langles D_1\rangles$. 
Then both triangles $\{\gamma, \gamma',\beta'\}$ and 
$\{\gamma, \gamma'' ,\beta''\}$ justify the admissibility
of the pair $(\gamma, \delta_1)$.}

\emph{Proof of the claim.} The condition (i) from the definition of admissibility 
in the section \ref{encoding-add} holds because both triangles $\{ \gamma, \gamma',\beta'\}$ and 
$\{\gamma, \gamma'' ,\beta''\}$ are marked.

Let us verify the condition (ii) from the definition of admissibility.
Since $D_1$ is a circle in $Q$, it is dijoint from the boundary $\partial Q{\eeq}C'_2\cup C'1_2$ and,
in particular, is disjoint from $C'_2$ and $C''_2$.
In addition, since $D_1$ is contained in $Q$ and $Q$ is disjoint from $C_1$,
the circle $D_1$ is also disjoint from $C_1$.
Since $\beta'$, $\beta''$ are the vertices of $\tts$ corresponding 
to the bounding pairs $\{C_1,C'_2\}$ , $\{C_1,C''_2\}$ respectively,
these disjointness properties of $D_1$ imply 
that $\delta_1$ is connected by edges of $\tts$ with both $\beta'$, $\beta''$.
This implies the condition (ii). 

Finally, let us verify the condition (iii).
By the choice of $D_1$ it is not isotopic in $Q$ to a circle disjoint from $C_0$.
As is well known, this implies that $D_1$ is not isotopic to such a circle in $S$ either.
It follows that $\delta_1$ is not connected by an edge of $\tts$ with $\gamma'$ or $\gamma''$.
This implies the condition (iii). This completes the proof of the claim. \esubproof
 
The above claim immediately implies that there is a move of
the type II replacing $(\gamma, \delta')$ by $(\gamma, \delta_1)$
and another move of the type II replacing $(\gamma, \delta_1)$ by 
$(\gamma, \delta'')$. Therefore, $(\gamma, \delta')$ is connected
to $(\gamma, \delta'')$ by two moves. This proves the lemma in Case 1. \eproof

\mytitle{Case 2. The circles\/ $C'_2$ and\/ $C''_2$ are contained in the same
component of $S\llsetmiss (C_0\cup C_1)$).} We may assume that this component is $S'$.
There are two subcases to consider depending on the genus of $S''$ (the part \emph{not}\/
containing $C'_2$ and\/ $C''_2$).

\mytitle{Subcase 2-A. The genus of\/ $S''$ is\/ $\geqs 2$.} Since 
the genus of\/ $S''$ is assumed to be\/ $\geqs 2$, there is
a non-trivial circle $E$ in $S''$ separating two boundary components $C_0$ and $C_1$ of $S''$.
Such a circle is not isotopic to $C_1$ in $S''$ (because it is non-trivial
in $S''$) and hence is not isotopic to $C_1$ in $S$. 
In addition, it is not isotopic to $C'_2$ in $S$ because $C'_2$ is contained in $S'$.
Finally, while the circle $E$ is separating in $S''$, it is non-separating in $S$,
because two sides of $E$ can be connected (by paths) in $S''$ to $C_0$ and $C_1$,
and $C_0$ and $C_1$ can be connected in $S'$.

Let $U_0$, $U_1$ be the components of $S''\llsetmiss E$ containing, respecively,
$C_0$ and $C_1$. Let $S_0{\eeq}U_0\cup C_0$ and $S_1{\eeq}U_1\cup C_1$. Clearly,
$S_0$ and $S_1$ are subsurfaces of $S$ having the unions $C_0\cup E$ and $C_1\cup E$
as, respectively, their boundaries. In addition, $C_0\cup E$ separates $S$ into subsurfaces
$S_1$ and $S'\cup S_0$, and $C_1\cup E$ separates $S$ into subsurfaces
$S_0$ and $S'\cup S_1$.

By the assumption, the union $C_0\cup C_1$ is separating $S$ (into $S'$ and $S''$).
Therefore, all three pairs $C_0\cup C_1$, $C_0\cup E$, and $C_1\cup E$, are bounding pairs
of circles, and the corresponding vertices $\gamma$, $\gamma'$ and $\beta$ of $\tts$
form a marked triangle.

Let $Q{\eeq}S'\cup S_0$. By exactly the same reasons as in Case 1, 
there is a non-trivial separating circle $D_1$ in $Q$
which is not isotopic in $Q$ to a circle disjoint from $C_0$.
Such a circle $D_1$ is not isotopic to a circle disjoint from $C_0$ also in $S$.
Since  $D_1$ is a circle in $Q$, it is disjoint from $C_1$ and $E$.
Let $\delta_1{\eeq}\langle D_1\rangle$. 
Obviously, the pair $(\gamma, \delta_1)$ encodes the same
isotopy class $\gamma_0{\eeq}\langles C_0\rangles$ as $(\gamma, \delta)$ and $(\gamma,\delta')$. 
Since the circle $D_1$ is disjoint form $C_1$ and $E$,
the vertex $\delta_1$ is connected with $\beta$ (recall that $\beta$ is the bounding pair
corresponding to the $C_1$, $E$) by an edge in $\tts$. It is clear now
that the triple $(\gamma,\gamma',\beta)$ justifies the admissibility of $(\gamma, \delta_1)$

Since $E$ is contained in $S''$ and both circles $C'_2$, $C''_2$ are contained in $S'$,
the already proved Case 1 of the lemma implies that
of the union $C_0\cup C_1$, we can use the result of the previous
paragraph and conclude that one can connect $(\gamma, \delta)$ with
$(\gamma, \delta_1)$ by two moves and also connect $(\gamma, \delta_1)$ with
$(\gamma, \delta')$ by two moves. It follows that one can connect
$(\gamma, \delta)$ with $(\gamma, \delta')$ by four moves. All these
moves are moves of the second type.  This proves the lemma in Subcase
2-A. \eproof

\mytitle{Subcase 2-B. The genus of\/ $S''$ is\/ $1$.} Since $S''$ has exactly two
boundary components, in this case $S''$ is a torus with two holes.
Since the surface $S$ is the union of its subsurfaces $S'$ and $S''$ intersecting
only along $\partial S'{\eeq}\partial S''{\eeq}C_0\cup C_1$, we may consider $S$ as the result of
gluening of $S''$ to $S'$ along boundary circles $C_0$, $C_1$. 
As is well known, when we start with a connected surface with boundary 
and glue to it another connected surface along $2$ boundary components,
the genus increases by $1+h$, where $h$ is the genus of the glued surface.
In our case this means that the genus of $S$ is equal to the genus of $S'$ plus $2$.
Since the genus of $S$ is assumed to be $\geqs 5$, this implies that the genus of $S'$ is $\geqs 3$. 

Since $\{C_0, C'_2\}$ is a bounding pair of circles, $C_0\cup C'_2$ divides $S$ into
two subsurfaces. Let us denote them by $P$, $Q$.
Since $C_1$ is disjoint from both circles $C_0$, $C'_2$, one of the
subsurfaces $P$, $Q$ contains $C_1$ and hence contains a neighborhood of $C_1$ in $S$.
In particular, it contains some points of the interior of $S''$. 
We may assume that this is the subsurface $Q$. 
Then $Q$ contains the whole subsurface $S''$ 
(because $\partial Q$ is disjoint from the interior of $S''$). 
This implies that $P$ is contained in $S'$.
In turn, this implies that $\{C_0, C'_2\}$ is a bounding pair of circles in $S'$
and hence $C'_2$ separates the boundary components $C_0$, $C_1$ of $S'$ in $S'$.
By the same argument, $C''_2$ also separates the boundary components $C_0$, $C_1$ of $S'$ in $S'$.

Since the genus of $S'$ is $\geqs 3$, we can apply Theorem \ref{special connectedness}
to $S'$ in the role of $R$. By this theorem, there is a sequence
\[
\seq{E}{C'_2}{C''_2}{n}
\]
of non-trivial circles in $S'$ separating the boundary components of $S'$
and such that $E_i$ and $E_{i+1}$ are isotopic in $S'$ to disjoint circles 
for all $i{\eeq}1,2,\ldots,n{\minus}1$. As usual, we may assume that the circles 
$E_i$ and $E_{i+1}$ are actually disjoint for all $i{\eeq}1,2,\ldots,n{\minus}1$.

For each $i{\eeq}1,2,\ldots,n$ let us denote by $R_i$ the subsurface of $S'$ bounded by $C_0$ and $E_i$. 
Since $E_i$ is nontrivial in $S'$, the genus of $R_i$ is $\geqs 1$.
This implies that the genus of $S''\cup R_i$ is $\geqs 2$.
Obviously, $C_0$ is a separating circle in $S''\cup R_i$ 
(it separates $S''\cup R_i$ into two parts $S''$ and $R_i$).
Since the genus of $S''\cup R_i$ is $\geqs 2$, this implies
that there is a non-separating circle $D_i$ in $S''\cup R_i$
which is not isotopic to a circle disjoint from $C_0$. 
Let $\delta_i{\eeq}\langles D_i\rangles$.  

For each $i{\eeq}1,2,\ldots,n$ the pairs $\{C_0,C_1\}$, $\{C_0,E_i\}$, and $\{C_1,E_i\}$
are bounding pairs of circles (recall that $E_i$ separates $C_0$ from $C_1$ in $S'$).
Recall that $\gamma$ is the vertex of $\tts$ corresponding to the bounding pair of circles $\{C_0,C_1\}$. 
Let $\gamma_i$, $\beta_i$ be the vertices of $\tts$ corresponding to the bounding pair of circles
$\{C_0,E_i\}$, $\{C_1,E_i\}$ respectively. 
Obviously, $\{\gamma,\gamma_i,\beta_i\}$ is a marked triangle in $\tts$.

\textsc{Claim.} \emph{For each\/ $i{\eeq}1,2,\ldots,n$ the pair\/ $(\gamma,\delta_i)$ is admissible.
In fact, the triple\/ $(\gamma,\gamma_i,\beta_i)$ justifies the admissibility of\/ $(\gamma,\delta_i)$.}

\emph{Proof of the claim.}\/ Since $\{\gamma,\gamma_i,\beta_i\}$ is 
a marked triangle for each\/ $i{\eeq}1,2,\ldots,n$, the condition (i)
from the section \ref{encoding-add} holds. 

For each $i{\eeq}1,2,\ldots,n$, the circle $D_i$ is a circle in $S''\cup R_i$
and hence is disjoint from $\partial (S''\cup R_i){\eeq}E_i\cup C_1$.
Therefore, $\delta_i$ is connected with $\beta_i$ by an edge of $\tts$
and the condition (ii) from the section \ref{encoding-add} holds. 

By the choice of the circle $D_i$, it is not isotopic to a circle disjoint from $C_0$ 
in $S''\cup R_i$, and hence also in $S$. This implies that $\delta_i$ is not connected
by an edge of $\tts$ with either $\gamma$ or $\gamma_i$. 
Therefore, the condition (iii) holds.
This completes the proof of the claim. \esubproof

\textsc{A sequence of moves.} Finally, we can construct a sequence of moves of type II
starting at $(\gamma,\delta')$ and ending at $(\gamma,\delta'')$.

Since $E_1{\eeq}C'_2$ and $E_n{\eeq}C''_2$, we have have $\gamma_1{\eeq}\gamma'$, $\gamma_n{\eeq}\gamma'$,
$\beta_1{\eeq}\beta'$, and $\beta_n{\eeq}\beta''$. Hence $(\gamma,\gamma_1,\beta_1){\eeq}(\gamma,\gamma',\beta')$
and $(\gamma,\gamma_n,\beta_n){\eeq}(\gamma,\gamma'',\beta'')$.
It follows that $(\gamma,\gamma_1,\beta_1)$ justifies the admissibilty of $(\gamma,\delta')$
in addition to $(\gamma, \delta_1)$. 
Hence $(\gamma,\delta')$ is connected with $(\gamma, \delta_1)$ by a move of type II.
Similarly, $(\gamma,\gamma_n,\beta_n)$ justifies the admissibilty of $(\gamma,\delta'')$
in addition to $(\gamma, \delta_n)$ and 
hence $(\gamma,\delta_n)$ is connected with $(\gamma, \delta'')$ by a move of type II.

For each $i{\eeq}1,2,\ldots,n{\minus}1$ the are two possible cases: the circle $E_{i+1}$ may be contained in $R_i$ or not. 
In the first case $D_{i+1}$ is a circle not only in $S''\cup R_{i+1}$, but also in $S''\cup R_i$.
In addition, $D_{i+1}$ is not isotopic to a circle disjoint from $C_0$.
Therefore, one can use $D_{i+1}$ in the role of $D_i$, and the above claim
implies that the triple $(\gamma,\gamma_i,\beta_i)$ justifies the admissibility of the pair
$(\gamma,\delta_{i+1})$ in addition to the pair $(\gamma,\delta_{i})$.
In particular, in this case the pairs $(\gamma,\delta_{i})$ and $(\gamma,\delta_{i+1})$
are connected by a move of type II.
Similarly, in the second case one can use $D_{i}$ in the role of $D_{i+1}$
the triple $(\gamma,\gamma_{i+1},\beta_{i+1})$ justifies the admissibility of the pair
$(\gamma,\delta_i)$ in addition to the pair $(\gamma,\delta_{i+1})$.
In particular, in this case the pairs $(\gamma,\delta_{i{\eeq}1})$ and $(\gamma,\delta_{i})$
are connected by a move of type II.

By combing the results of the last two paragraphs, we conclude that there is
a sequence of moves of type II starting at $(\gamma,\delta')$ and ending at $(\gamma,\delta'')$.
This completes the proof in the case 2-B, and hence of the lemma \eproof

\mytitle{Remark.} The assumption that the genus of $S$ is $\geqs 5$ is used at
only one place in the above proof. Namely, it is used only in the first paragraph
of the proof in the Subcase 2-B in order to ensure that the genus of $S'$ is $\geqs 3$.
It is fairly reasonable to expect that the proof can be extended at
least to the case of surfaces of genus $4$ at the cost of introducing some 
additional types of moves.

\mypar{Theorem \textup{(}Connecting encoding pairs by moves\textup{)}.}{changing gamma and delta} \emph{Suppose
that the genus of $S$ is $\geqs 5$.
Suppose that a non-separating circle can be encoded by two (different) admissible pairs
$(\gamma, \delta)$ and $(\gamma', \delta')$. Then there is a sequence 
of moves of types I and II, starting at $(\gamma, \delta)$ and ending at $(\gamma', \delta')$.}

\emph{Proof.}\/ Let $C_0$ be a non-separating circle such that $\gamma_0{\eeq}\langles C_0\rangles$ 
is the isotopy class encoded by the both pairs. 
Then for some no-separating circles $C_1$, $C'_1$ the vertices $\gamma$, $\gamma'$ of $\tts$ 
correspond to bounding pairs of circles  $\{C_0, C_1\}$ and $\{C_0, C'_1\}$  respectively.

Let $R$ be the result of cutting $S$ along $C_0$. 
Clearly, $R$ is surface with two boundary components, 
and its genus is equal to the genus of $S$ minus $1$. In particular,
the genus of $R$ is $\geqs 4$.
By Theorem \ref{special connectedness} there is a sequence
\[
\seq{E}{C_1}{C'_1}{n}
\] 
of nontrivial circless in $R$ separating its boundary components and 
and such that $E_i$ and $E_{i+1}$ are isotopic in $R$ to disjoint circles 
for all $i{\eeq}1,2,\ldots,n{\minus}1$. As usual, we may assume that the circles 
$E_i$ and $E_{i+1}$ are actually disjoint for all $i{\eeq}1,2,\ldots,n{\minus}1$.
Let $\gamma_i$ be the vertex of $\tts$ corresponding to the bounding pair of circles $C_0$, $E_i$. 
In particular, $\gamma_1{\eeq}\gamma$ and $\gamma_n{\eeq}\gamma'$.
 
\textsc{Claim.} \emph{There are non-trivial circles\/ $D_1$, $D_2$, \ldots, $D_n$, such that
for each\/ $i{\eeq}1,2,\ldots,\mbox{ or }n{\minus}1$:}
\vspace{-3ex}
\begin{itemize}
\item[(i)] \emph{the pair\/ $(\gamma_i,\delta_i)$, where\/ $\delta_i{\eeq}\langles D_i\rangles$, encodes\/ $\gamma_0$\,;} 
\vspace{-1ex}
\item[(ii)] \emph{$(\gamma_i, \delta_i)$ is connected with\/ $(\gamma_{i+1}, \delta_i)$ by a move of type I\/\/\/;}
\vspace{-1ex}
\item[(iii)] \emph{$(\gamma_{i+1}, \delta_i)$ is connected with\/ $(\gamma_{i+1}, \delta_{i+1})$ 
by a sequence of moves of type II\/\/\/;}
\vspace{-1ex}
\item[(iv)] \emph{$(\gamma, \delta)$ is connected with\/ $(\gamma, \delta_1){\eeq}(\gamma_1, \delta_1)$ by a
sequence of moves of type II\/\/\/;}
\vspace{-1ex}
\item[(v)] \emph{$(\gamma_n, \delta_n){\eeq}(\gamma', \delta_n)$ is connected with\/ $(\gamma', \delta')$ by a
sequence of moves of type II\/.}
\end{itemize}
\vspace{-3.5ex}
\vspace{\parskip}

\emph{Proof of the claim.}\/ Let $i{\eeq}1,2,\ldots,\mbox{ or }n{\minus}1$.
 
Since both $E_i$ and $E_{i+1}$ separate the boundary components of $R$,
the union $E_i\cup E_{i+1}$ divides $R$ into $3$ parts.
One of these parts is a subsurface of $R$ with boundary $E_i\cup E_{i+1}$.
We will denote this subsurface by $P_i$.
Two other parts are the components of the closure of $R\llsetmiss P_i$ in $R$.
Clearly, these components are subsurfaces of $R$ having two boundary components each.
One of these subsurfaces has $E_i$ as a boundary component, the other has $E_{i+1}$.
We will denote these subsurfaces by $V_i$ and $W_i$ 
(it does not matter which one is denoted by $V_i$ and which one by $W_i$).
Since $E_i$ and $E_{i+1}$ are non-trivial circles in $R$, neither of subsurfaces $V_i$, $W_i$ is an annulus.
Since each of surfaces $V_i$, $W_i$ has two boundary components, this implies
that the genus of surfaces $V_i$, $W_i$ is $\geqs 1$.

Clearly, $P_i$, $V_i$, and $W_i$ are also subsurfaces of $S$. 
While $V_i$ and $W_i$ are dijoint as subsurfaces of $R$, 
their intersection as subsurfaces of $S$ is equal to $C_0$. 
Let $Q_i$ be the union of $V_i$ and $W_i$ considered as subsurfaces of $S$.  
It is a subsurface of $S$ with the boundary $\partial Q_i$ equal 
to the boundary $\partial P_i {\eeq}E_i\cup E_{i+1}$.
Since the genus of surfaces $V_i$, $W_i$ is $\geqs 1$, the genus of $Q_i$ is $\geqs 2$.
Obviously, $C_0\subset Q_i$ and hence $C_0$ is a circle in $Q_i$.
Clearly, the set difference $Q_i\setmiss C_0$ consists of two components.
Namely, the components of $Q_i\setmiss C_0$ are  $V_i\setmiss C_0$ and $W_i\setmiss C_0$.
Since, as we saw, neither $V_i$, nor $W_i$ is an annulus,
the circle $C_0$ is a non-trivial separating circle in $Q_i$ 
(the circle $C_0$ does not bound a disc in $S$ and, hence, in $Q_i$).
Since the genus of $Q_i$ is $\geqs 2$, this implies that $Q_i$
contains a circle not isotopic to a circle disjoint from $C_0$.
Let us choose such a circle as the promised circle $D_i$. 
As in the statement of the claim, we will denote by $\delta_i$ 
the isotopy class $\langles D_i\rangles$.

Since the circle $D_i$ is contained in $Q_i$, it is disjoint from $C_1$.
In addition, $D_i$ is not isotopic to a circle disjoint from $C_0$.
This implies that $(\gamma_i,\delta_i)$ encodes $\gamma_0$, and hence imples (i).
 
Let $\beta_i$ be the vertex of $\tts$ corresponding to the bounding pair of circles $E_i$, $E_{i+1}$.
Clearly, $\{\gamma_i,\gamma_{i+1},\beta_i\}$ is a marked triangle.
Moreover, the triple $(\gamma_i,\gamma_{i+1},\beta_i)$ justifies the admissibility of both pairs
$(\gamma_i, \delta_i)$ and $(\gamma_{i+1}, \delta_i)$.
Therefore, $(\gamma_i, \delta_i)$ is connected with $(\gamma_{i+1}, \delta_i)$ by a move of type I.
This proves (ii).

Finally, Lemma \ref{changing delta} implies (iii), (iv), and (v) 
because in these statements only the second element of the encoding pair is changed.
This completes the proof of the claim. \esubproof

By combining the moves of type I from the statement (ii) of the claim
with the sequences of moves of type II from the statements (iii), (iv), and (v) of the claim,
we get a sequence of moves connecting $(\gamma, \delta)$ with $(\gamma', \delta')$.
This completes the proof of the lemma. \eproof

\mypar{Corollary.}{encoding} \emph{Suppose that the genus of\/ $S$ is\/ $\geqs 5$.
Two admissible pairs encode the same isotopy class of non-separating curve if 
and only if they can be connected by a sequence of moves of types I and II.}

\emph{Proof.} The corollary immediately follows from 
Lemmas \ref{not changed} and \ref{changing gamma and delta}.\eproof

\mysection{Encodings, edges of $\ccs$, and automorphisms}{en-edges-aut}

From now on, for the isotopy classes of non-separating circles we will consider 
only encodings by admissible pairs, and will call them simply \emph{encodings}.

Let $D$ be a separating circle in $S$, and let $\delta{\eeq}\llvv{D}$.
The isotopy class $\delta$ is a vertex of $\ccs$ and at the same time is a vertex of $\tts$.
It is convenient define an \emph{encoding of}\/ $\delta$ considered as a vertex of $\ccs$
as the same isotopy class $\delta$, but considered as a vertex of $\tts$.
In contrast with non-separating vertices, separating circles have only one encoding.
With this definition, we have the notion of an encoding for the isotopy
classes of all non-trivial circles in $S$, or, what is the same, 
for all vertices of $\ccs$.

The goal of this section is to find properties of encodings of vertices of $\ccs$
which are equivalent to the property of being connected in $\ccs$, and which
can be stated entirely in terms of the Torelli building $\tts$. As an application,
we will prove Theorem \ref{automorphisms of geometry} at the end of this section.
See the subsection \ref{proof-aut-tg}.

\mypar{Recognizing circles of genus $1$.}{genus-1} Let $\delta$ be a vertex of $\ccs$,
i.e. $\delta {\eeq} \vv{D}$ for a non-trivial circle $D$ in $S$. 
For the problem of deciding if $\delta$ is connected with another vertex $\gamma$ of $\ccs$ 
by and edge of $\ccs$ in terms of encodings of $\delta$, $\gamma$, 
the case of genus $1$ circles $D$ presents substantial 
additional difficulties (in comparison with other circles in $S$). 
The main reason for this is the fact that a torus with one hole 
contains neither non-trivial separating circles, nor bounding pairs of circles. 

We will need to deal with such vertices $\delta {\eeq} \vv{D}$ separately.
In order to do this, we, first of all, need to be able to recognize
such vertices in terms of properties their encodings in the Torelli building $\tts$. 
In more details, we need to find properties of encodings of a vertex $\delta {\eeq} \vv{D}$ of $\ccs$
which are equivalent to $D$ being a genus $1$ circle, and which are stated entirely in terms of the Torelli building $\tts$. 
This is the main goal of this subsection, achieved in Lemma \ref{links} below.

The first step in achieving this goal is trivial: encodings of separating and non-separating
vertices are of different nature. Therefore, $D$ can be assumed to be a separating circle from
the very beginning. Still, we will use this this assumption only when it is needed.

\myitpar{The graph $\link\sep(\delta)$.} Our main tool will be 
the graph $\link\sep(\delta)$ which we will define now.
Let $L\sep(\delta)$ be the set of separating vertices of $\tts$
different from $\delta$ and connected by an edge of $\tts$ with $\delta$.
Let us turn $L\sep(\delta)$ into a graph $\link(\delta)$ 
by connecting two vertices in $L\sep(\delta)$ by
an edge if and only they are \emph{not}\/ connected by an edge in $\tts$.
The graph $\link\sep(\delta)$ is called the \emph{separating link}\/ of $\delta$ (in $\tts$).

In order to deal with graphs $\link\sep(\delta)$, we need auxillary graphs $G\sep(R)$ 
defined for a compact connected surfaces $R$ with exactly one boundary component. 
Let $\mbox{V}\sep(R)$ be the set of isotopy classes of separating circles on $R$. 
Clearly, $\mbox{V}\sep(R)$ is a subset of the set of vertices of $\ccc(R)$. 
Let us define $G\sep(R)$ as the graph having $\mbox{V}\sep(R)$ as its set of vertices 
and having two different vertices connected by an edge if and only if these vertices are 
\emph{not connected by an edge}\/ of $\ccc(R)$.

\mylemma{Observation.}{sep-with-1} \emph{Let $Q$ be a compact connected surface with
$1$ boundary component, and let $D$ be a separating circle in $Q$. Then one of the subsurfaces
into which $D$ divides $Q$ has $D\cup\partial Q$ as its boundary, and the other is
disjoint from $\partial Q$ and its boundary is $D$.}

\proof Since $\partial Q$ is connected, one of these subsurfaces contains $\partial Q$.
Clearly, the boundary of this subsurface is $D\cup\partial Q$. 
It follows that the other subsurface is disjoint from $\partial Q$ 
and has $D$ as its boundary. \eproof

\mylemma{Lemma.}{sep-r-q} \emph{Suppose that the surface $R$ is a subsurface of 
another surface $Q$ with $\leqs 1$ boundary component.
Let $D$ be a circle contained in $R\llsetmiss\partial R$. 
Then the circle $D$ is a separating circle in $R$
if and only if it is a separating circle in $Q$.}

\proof Suppose that $D$ is a separating circle in $R$. 
By Observation \ref{sep-with-1}, $D$ divides $R$ into two subsurfaces, and
one of these subsurfaces has $D\cup\partial R$ as its boundary, 
while the other is disjoint from $\partial R$ and has $D$ as its boundary. 
The second subsurface is also a subsurface of $Q$.
It follows that $D$ is a separating circle in $Q$ . 
This proves the \emph{``only if''} direction of the lemma. 

By the same Observation \ref{sep-with-1}, if $D$ is a separating circle in $Q$, then $D$ divides $Q$ 
into two subsurfaces such that one of them has $D\cup\partial Q$ as its boundary, 
while the other one has $D$ as its boundary.
Clearly, if $P$ is a subsurface of $Q$ disjoint from $\partial Q$ with $\partial P {\eeq} D \subset R$,
then $P$ is contained in $R$. It follows that if $D$ is a separating circle in $Q$, 
the $D$ is the boundary of a subsurface of $R$ and hence $D$ is a separating circle in $R$.
This proves the \emph{``if''} direction of the lemma, and hence completes the proof. \eproof

\mylemma{Lemma.}{connecteness-of-dual} \emph{Let\/ $R$ be a compact connected surface with\/ $1$ boundary component.
If\/ $\genus(R)\geqs 2$, then\/ $\mbox{\rm V}\sep(R)\neq\varnothing$ and
the graph\/ $G\sep(R)$ is connected.}

\proof It is much more easy to connect two vertices in 
$G\sep(R)$ than to connect them in $\ttg(R)$.
The experts may skip the rest of the proof. 

The non-emptiness statement is trivial. Note that the assumption $\genus(R)\geqs 2$ is exactly
what is needed for it to be true (for surfaces with one boundary component).

In order to prove the connectedness statement of the lemma, 
it is sufficient to prove that the isotopy classes $\delta{\eeq}\llvv{D}$, $\delta'{\eeq}\llvv{D'}$ 
of two arbitrary non-trivial separating circles $D$, $D'$ in $R$ are connected by
a sequence of edges of $G\sep(R)$.

\mytitle{Case 1. The circles\/ $D$, $D'$ are not isotopic to disjoint circles.}\/ 
Then their isotopy classes $\delta{\eeq}\llvv{D}$, $\delta'{\eeq}\llvv{D'}$ 
are not connected by an edge of $\ccc(R)$, 
and hence are connected by an edge of $G\sep(R)$. 
This completes the proof in Case 1. \eproof

\mytitle{Case 2. The circles\/ $D$, $D'$ are isotopic to disjoint circles.}\/ 
In this case we may assume that the circles $D$, $D'$ are themselves disjoint. 
Since $D$ is a separating circle in $R$, it divides $R$ into two subsurfaces 
such that one of them has $D\cup\partial R$ as its boundary, 
while the other has $D$ as its boundary and is disjoint from $\partial R$.
Let $Q_1$ be the first subsurface, and let $Q$ be the second, so, in particular, $\partial Q {\eeq} D$. 
Since $D'$ is disjoint from $\partial Q{\eeq}D$, either $D'\subset Q$ or $D'$ is disjoint from $Q$. 

\mytitle{Subcase 2-A. $D'\subset Q$.}\/ Since $D'$ is a separating circle in $R$, Lemma \ref{sep-r-q}
implies that $D'$ is also a separating circle in $Q$. This implies that $D'$ 
divides $Q$ into two subsurfaces having $D$ as the common component of their boundaries.
One of them has $\partial Q$ as a boundary component, the other is disjoint from $\partial Q$.
Let $Q_3$ be the first subsurface, and let $Q_2$ be the second. 
In particular, $\partial Q_3 {\eeq} D'$, and $\partial Q_2 {\eeq} D'\cup \partial Q {\eeq} D'\cup D$.

It follows that $R$ is the union of $3$ subsurfaces $Q_3$, $Q_2$, $Q_1$ with boundaries
$\partial Q_3{\eeq}D'$, $\partial Q_2{\eeq}D'\cup D$, and $\partial Q_1 {\eeq} D\cup\partial R$.
Moreover, $Q_3\cap Q_2 {\eeq} D'$, $Q_2\cup Q_1 {\eeq} D$, and $Q_3\cap Q_1 {\eeq} \varnothing$.
It follows that the configuration of surfaces $R$, $Q_3$, $Q_2$, $Q_1$ and circles
$D$, $D'$ up to diffeomorphisms depends only on the genera 
$\genus(Q_3)$, $\genus(Q_2)$, and $\genus(Q_1)$.
Therefore, we may assume that we are dealing with a standard model configuration.
Then it is very easy to exhibit a separating circle $A$ in $R$ such that $A$ is
not isotopic to a circle disjoint from either $D$ or $D'$. 
We leave this task to the reader.
If $A$ is such a circle, then $\alpha{\eeq}\llvv{A}$ is not connected by edges of $\ccc(R)$ 
with $\delta{\eeq}\llvv{D}$ and $\delta'{\eeq}\llvv{D'}$. By the definition, this means that
$\alpha$ is connected with both $\delta$ and $\delta'$ by edges of $G\sep(R)$.
Hence, $\delta$ is connected with $\delta'$ by a sequence consisting of these two edges. 
This completes the proof in Subcase 2-A. \eproof

\mytitle{Remark.} Since circles $D$, $D'$ are non-trivial, in Subcase 2-A 
either $D$ is isotopic to $D'$,
or each of the subsurfaces $Q_3$, $Q_2$, $Q_1$ has genus $\geqs 1$.
Therefore, in Subcase 2-A either $\delta{\eeq}\delta'$ (and then there is
nothing to prove), or $\genus(R)\geqs 3$.

\mytitle{Subcase 2-B. $D'$ is disjoint from $Q$.}\/
Since $D'$ is a separating circle in $R$, $D'$ divides $R$ into two subsurfaces 
such that one of them has $D'\cup\partial R$ as its boundary, 
while the other has $D'$ as its boundary and is disjoint from $\partial R$.
Let $Q'_1$ be the first subsurface, and let $Q'$ be the second, so, in particular, $\partial Q' {\eeq} D'$.
Since $D$ is disjoint from $\partial Q'{\eeq}D'$, either $D\subset Q'$ or $D$ is disjoint from $Q'$.
If $D\subset Q'$, then the arguments of the Subcase 2-A apply, with the roles of $D$ and $D'$ interchanged. 
Hence we may assume that $D$ is disjoint from $Q'$.

\textsc{Claim.} \emph{If\/ $D'$ is disjoint from\/ $Q$, and\/ $D$ is disjoint from\/ $Q'$,
then\/ $Q$ is disjoint from\/ $Q'$.}

\emph{Proof of the claim.} Consider the union $P {\eeq} Q\cup Q'$. 
Since the boundaries $D {\eeq} \partial Q$ and $D' {\eeq} \partial Q'$ are 
disjoint, $P$ is subsurface of $R$ with the boundary $D\cup D'$. 
If $P$ is connected, then either $P\subset Q$, or $P\subset Q_1$. 
Since $D'\subset\partial P\subset P$,
the inclusion $P\subset Q$ contradicts to the assumption
that $D'$ is disjoint form $Q$.
If $P\subset Q_1$, then $Q\subset P \subset Q_1$, contradicting to the
fact that $Q\cap Q_1 {\eeq} D$.
It follows that $P$ is not connected. Since $Q$ and $Q'$ are connected,
this is possible only if $Q$ and $Q'$ are disjoint.
This proves our claim. \esubproof

If $Q$ is disjoint from $Q'$, the the closure of $R\llsetmiss (Q\cup Q')$
is a suburface of $R$ with the boundary $D\cup D'\cup\partial R$. 
Let us denote this subsurface by $Q_0$. 
Then $R$ is the union of $3$ subsurfaces $Q$, $Q'$, and $Q_0$ with boundaries
$\partial Q{\eeq}D$, $\partial Q'{\eeq}D'$, and $\partial Q_0 {\eeq} D\cup D'\cup\partial R$.
Moreover, $Q\cap Q' {\eeq} \varnothing$, $Q\cup Q_0 {\eeq} D$, and $Q'\cap Q_0 {\eeq} D'$.
It follows that the configuration of surfaces $R$, $Q$, $Q'$, $Q_0$ and circles
$D$, $D'$ up to diffeomorphisms depends only on the genera 
$\genus(Q)$, $\genus(Q')$, and $\genus(Q_0)$.
Therefore, we may assume that we are dealing with a standard model configuration.
Then it is very easy to exhibit a separating circle $A$ in $R$ such that $A$ is
not isotopic to a circle disjoint from either $D$ or $D'$. 
We leave this task to the reader. Given such a circle $A$, we can finish the
proof in Subcase 2-B in the same manner as we did in the Subcase 2-A.
This completes the proof in Subcase 2-B and, hence, the proof of the lemma. \eproof

\mylemma{Lemma.}{links} \emph{Suppose that $\gen{\eeq}\genus(S)\geqs 3$.
Let $D$ be a separating circle in $S$, 
and let $\delta{\eeq}\langles D\rangles$. 
Then the graph\/ $\link\sep(\delta)$ is connected if and only if\/ 
$D$ bounds in\/ $S$ a torus with one hole.}

\proof Let $S_1$ and $S_2$ be the two parts into which $D$ separates $S$. 
Both $S_1$ and and $S_2$ are subsurfaces of $S$ with boundary $D$.
Every vertex of $\link\sep(\delta)$ can be represented by a separating circle $A$ disjoint from $D$. 
Obviously, every such a circle $A$ is contained either in $S_0$ or in $S_1$.
Since $\vv{A}\neq\delta$, the circle $A$ does not bound an annulus in $S$ together with $D$
and hence if $A\subset S_i$, then $A$ is a non-trivial circle in $S_i$ (where $i{\eeq}0$ or $1$).

It follows that the set of vertices of $\link\sep(\delta)$ is equal to the union $L_1\cup L_2$
of sets $L_1$, $L_2$, where for each $i{\eeq}1,2$ 
the set $L_i$ consists of the vertices which can be represented
by separating non-trivial circles in $S_i\setmiss D$. 
Clearly, such a representative of any vertex in $L_1$
is disjoint from such a representative of any vertex in $L_2$.
Therefore, every element of $L_1$ is connected by an edge of $\ccs$ with every element of $L_2$. 
By the defintion of $\link\sep(\delta)$, this means that vertices 
from $L_1$ are not connected with vertices of $L_2$ by edges of $\link\sep(\delta)$.
Therefore, if both sets $L_1$, $L_2$ are non-empty, then 
the graph $\link\sep(\delta)$ has at least $2$ connected components.

Together with Lemma \ref{connecteness-of-dual} this implies that if $G(S)$ is connected, 
then either $\genus(S_1)\leqs 1$ or $\genus(S_2)\leqs 1$.
We may assume that $\genus(S_1)\leqs 1$. 
Since $S_1$ is a connected surface with exactly $1$ boundary component 
and is not a disc (because $\partial S_1{\eeq}D$ is a non-trivial circle in $S$), 
the classification of surfaces implies that in this case $S_1$ is a torus with one hole. 
This proves the \emph{``only if''}\/ direction of the lemma. 

In order to prove the \emph{``if''}\/ direction of the lemma, 
we may assume that $S_1$ is a torus with one hole. 
Then $L_1{\eeq}\varnothing$, the set of vertices of the graph $\link\sep(\delta)$ is equal to $L_2{\eeq}\mbox{V}\sep(S_2)$, 
and hence $\link\sep(\delta){\eeq}G\sep(S_2)$.
Since $S_1$ is a torus with one hole, $\genus(S_2){\eeq}\genus(S){\minus}1{\eeq}\gen{\minus}1\geqs 2$.
Therefore, Lemma \ref{connecteness-of-dual} applies to $S_2$ 
and implies that $\link\sep(\delta){\eeq}G\sep(S_2)$ is connected.
This completes the proof of the \emph{``if''}\/ direction of the lemma, 
and hence completes the proof. \eproof

\mypar{Circles of genus $1$ and non-separating circles.}{genus-one-non-sep} In this subsection we assume that $\genus(S)\geqs 3$.
The main results of this subsection, namely Lemmas \ref{c-in-t-encode} and \ref{c-out-t-encode}, 
require stronger assumption $\genus(S)\geqs 4$.
 
Let $D_0$ be a genus $1$ circle in $S$, and let $\delta_0 {\eeq} \llvv{D_0}$. 
Since $\genus(S)\geqs 3$, only one of the two parts into which 
$D_0$ divides $S$ can be a torus with one hole.
Let $S_0$ be this part, and let $S_1$ be the other. 
Our next goal is to find for $i{\eeq}0,1$ properties of encodings of isotopy classes 
$\gamma {\eeq}\llvv{C}$ of non-separating circles $C$ in $S$
which are equivalent to $C$ being isotopic to a circle contained in $S_i$,
and which can be stated entirely in terms of the 
Torelli building $\tts$ and its vertex $\delta$.
Such conditions are provided by Lemmas \ref{c-in-t-encode} and \ref{c-out-t-encode} below.

\mylemma{Lemma.}{same-part} \emph{Suppose that $D$ is a non-trivial separating circle in $S$.
If\/ $C_0$, $C_1$ is a bounding pair of circles disjoint from\/ $D$,
then both circles\/ $C_0$, $C_1$ are contained in the same component of\/ $S\llsetmiss D$.}

\proof Let $S_0$ and $S_1$ be the closures of two components of $S\llsetmiss D$.
Clearly, both $S_0$ and $S_1$ are subsurfaces of $S$ having $D$ as their common boundary.
Being connected, each of the circles $C_0$, $C_1$ is entirely contained either in $S_0$ or in $S_1$.

Suppose that $C_0\subset S_0$. If $C_0$ separates $S_0$, 
then $C_0$ bounds in $S_0$ a subsurface not containing $D$, 
and then $C_0$ bounds the same subsurface in the whole surface $S$.
This contradicts to the fact that $C_0$ is non-separating in $S$, 
being a circle from a bounding pair of circles in $S$.
It follows that $C_0$ is a non-separaring circle in $S_0$. 
By the same argument $C_1$ is a non-separating circle either in $S_0$ or in $S_1$,
depending on which of these subsurfaces contains $C_1$.

Recall that $C_0\subset S_0$ by our assumption. If $C_1\subset S_1$,
then the union $C_0\cup C_1$ is does not separates $S$
because circles $C_0$, $C_1$ are non-separating in $S_0$, $S_1$ respecitively,
and subsurfaces $S_0$, $S_1$ are connected by their common boundary.
This contradicts to the fact that $C_0$, $C_1$ is a bounding pair of circles.
It follows that $C_1\subset S_0$.

Therefore, if $C_0\subset S_0$, then also $C_1\subset S_0$ 
Similarly, if $C_0\subset S_1$, then $C_1\subset S_1$. 
The lemma follows. \eproof

\mytitle{Remark.} The last lemma is valid without any restrictions on $\gen{\eeq}\genus(S)$,
but it is non-vacuous only if $\gen\geqs 2$.

\mylemma{Lemma.}{circle-in-torus} \emph{Let\sss $C_0$ be a non-separating circle in\sss $S$.\sss
The circle\sss\sss $C_0$ is isotopic to a circle contained in\sss $S_0$ if and only if
for each separating circle\sss $D$\sss in\sss $S_1$, the circle\sss $C_0$ is isotopic to a circle disjoint from\sss $D$.}

\proof The \emph{``only if''} direction of the lemma is trivial. 
In order to prove the \emph{``if''} direction of the lemma, 
we will use Thurston's theory of surfaces.
In fact, we will prove a much stronger result.

Recall that a collection of circles $C_1\/, C_2\/,\ldots\/,C_n$ contained in $S_1$ is said to be
\emph{filling the surface}\sss $S_1$ if for every circle $C$ in $S$ the following two conditions are equivalent.
\vspace{-2.7ex}
\begin{itemize}
\item[(i)] \emph{For each for\sss $i{\eeq}1,2,\ldots,n$\sss the circle\sss $S$\sss is isotopic a circle disjoint from\sss $S_1$.}
\vspace{-1.5ex}
\item[(ii)] \emph{The circle\sss $C$\sss is isotopic to a circle disjoint from\sss $S_1$.}
\end{itemize}
\vspace{-2.7ex}
Note that, obviously, (ii) always implies (i).
Thurston's theory provides a powerful tool for constructing filling collections of circles.
Actually, every connected surface can be filled by just two circles.

Our assumption $\genus(S)\geqs 3$ implies that $\genus(S_1)\geqs 2$, 
and hence $S_1$ contains non-trivial separating circles 
(see, for example, Lemma \ref{connecteness-of-dual}). 
Let $E$ be such a circle. 
Let $f{\in}\Mod(S_1)$ be a pseudo-Anosov element,
and let $F\colon S_1{\tto} S_1$ be any diffeomorphism representing $f$. 
Then for any sufficiently large natural numbers $n$, $m$
the circle $F^n(E)$, $F^{-m}(E)$ fill $S_1$ 
(the fact that $D$ is separating is of no importance here).
It follows that any circle in $S$ isotopic to a circle disjoint from $F^n(E)$,
and isotopic to a, perhaps, another circle disjoint from $F^{-m}(E)$
is isotopic to a circle disjoint from $S_1$ and hence contained in $S_0$.
This completes the proof of the \emph{``if''} direction of the lemma,
and hence proves the lemma. \eproof

\mylemma{Corollary}{c-in-t} \emph{Let\/ $C_0$ be a non-separating circle in\/ $S$, 
and let $\gamma_0{\eeq}\vv{C_0}$. The circle\/ $C_0$ in\/ $S$ is isotopic
to a circle contained in\/ $S_0$ if and only if for every separating circle\/ $D$ in\/ $S_1$
the vertex $\gamma_0{\eeq}\vv{C_0}$ of\/ $\ccs$ is connected by an edge of\/ $\ccs$
\emph{(}\hspace{0.05em}or, what is the same, of\/ $\tts$\hspace{0.05em}\emph{)} 
with\/ $\delta{\eeq}\llvv{D}$.} \eproof

\mytitle{Remark.} Unfortunately, the authors are not aware of any reference for 
the used above corollary of Thurston's theory not requiring learning the whole theory.
For whose who know the whole theory (with proofs), this corollary is obvious. 
The idea behind it is that the images of $\vv{F^{n}(E)}$ and $\vv{F^{-m}(E)}$ in 
the projective space of measured foliations approach, respectively, 
the stable and unstable foliations of $f$ when $n,m{\tto}\infty$.

\mylemma{Lemma.}{c-in-t-encode} \emph{Suppose that $\genus(S)\geqs 4$. 
Let\/ $C_0$ be a non-separating circle in\/ $S$,  and let\/ $\gamma_0{\eeq}\vv{C_0}$. 
The circle\/ $C_0$ is isotopic to a circle contained in\/ $S_0$ 
if and only if the following condition holds.}

\vspace{-1ex}
\hangindent=1em \hangafter=0
$({\mathcal S_0})$ \emph{For every non-trivial separating circle\/ $D$ in\/ $S_1$,
there exists a admissible encoding\/ $(\gamma,\varepsilon)$ of\/ $\gamma_0{\eeq}\vv{C_0}$ 
such that\/ $\gamma$ is connected with\/ $\delta{\eeq}\llvv{D}$ by an edge of\/ $\tts$.}

\proof Suppose that $C_0$ is contained in $S_0$. 
Let $D$ be a non-trivial separating circle in $S_1$.
Since $D$ is separating, $D$ divides $S_1$ into two parts.
One of these parts contains $D_0{\eeq}\partial S_1$.
Let us denote this part by $R$. 
Then $R$ is a subsurface of $S_1$ with boundary $\partial R {\eeq} D_0\cup D$.
Since $D$ is non-trivial, $D$ is not isotopic in $S$ to $D_0{\eeq}\partial S_1$.
Therefore, $R$ is not an annulus and hence $\genus(R)\,{\geq}\, 1$. 
This implies that the genus of the subsurface $R_0{\eeq}S_0\cup R$ is ${\geq}\, 2$.
Obviously, $C_0$ is a non-separating circle in $R_0$.
Since $\genus(R_0){\geq} 2$, there is a bounding pair of circles in $R_0$
having $C_0$ as one of its components.
Let $C_1$ be the other circle from this bounding pair of circles, and let
$\gamma{\eeq}\{\vv{C_0},\vv{C_1}\}$.
Since $\genus(S){\geq} 3$, there is a separating circle $E$ in $S$ such
that the pair $(\gamma,\varepsilon)$, where $\varepsilon{\eeq}\llvv{E}$, is an encoding of $\gamma_0$. 
Since $\genus(S)\geqs 4$, we can choose such a circle $E$ that this encoding is admissible.
Since $C_0\cup C_1$ is contained in $R_0$ and hence disjoint from $D{\eeq}\partial R_0$,
the vertices $\gamma$ and $\delta$ of $\tts$ are connected by an edge of $\tts$.
This proves the \emph{``only if''} part of the lemma.

In order to prove the \emph{``if''} part, suppose that $\gamma$ is connected with
$\delta$ by an edge of $\tts$ for an encoding $(\gamma,\varepsilon)$ of $\gamma_0$.
For every such an encoding the vertex $\gamma$ has the form 
$\gamma{\eeq}\{\vv{C_0},\vv{C_1}\}$ for some non-separating circle $C_1$,
and since $\gamma$ is connected by an edge with $\delta$, both circles $C_0$, $C_1$
are isotopic to circles disjoint from $D$. In particular, $C_0$ is isotopic to a
circle disjoint from $D$, and hence $\gamma_0$ is connected by an edge of $\ccs$ with $\delta$.

Therefore, $({\mathcal S_0})$ implies that for every separating circle $\delta$ in $S_1$ the vertex
$\gamma_0$ of $\ccs$ is connected with $\delta$ by an edge of $\ccs$.
Now Corrolary \ref{c-in-t} implies that $C_0$ is isotopic to a circle in $S_0$.
This proves the \emph{``if''} part of the lemma, and hence completes the proof. \eproof

\mylemma{Lemma.}{c-out-t-encode} \emph{Suppose that $\genus(S)\geqs 4$. 
Let\/ $C_0$ be a non-separating circle in\/ $S$,  and let\/ $\gamma_0{\eeq}\vv{C_0}$. 
The circle\/ $C_0$ is isotopic to a circle contained in\/ $S_1$ 
if and only if the following condition holds.}

\vspace{-1ex}
\hangindent=1em \hangafter=0
$({\mathcal S_1})$ \emph{There exists an admissible encoding\/ $(\gamma,\varepsilon$) of\/ $\gamma_0$ 
such that the vertices\/ $\delta_0$ and\/ $\gamma$ of\/ $\tts$ are connected by an edge of\/ $\tts$.}

\proof Note that $\genus(S_1){\eeq}\genus(S){-}1 \geqs 3$. 

Suppose that $C_0$ is isotopic to a circle in $S_1$. Then we can
assume that the circle $C_0$ itself is contained in $S_1$.
Since the genus $\genus(S_1)\geqs 3$, there exist a non-separating circle $C_1$
and a separating circle $E$ in $S_1$ such that $(\gamma,\varepsilon)$ is an encoding of $\gamma_0$,
where $\gamma{\eeq}\{\vv{C_0},\vv{C_1}\}$ and $\varepsilon{\eeq}\vv{E}$.
Since $\genus(S)\geqs 4$, the circles $C_1$ and $E$ can be chosen in a such a way that 
this encoding is admissible. Since $C_0$, $C_1$ are circles in $S_1$, 
they are disjoint from $D_0{\eeq}\partial S_1$, and hence $\gamma$ is connected by an edge of $\tts$ with $\delta_0$.
This proves the \emph{``only if''} direction of the lemma.

In order to prove the \emph{``if''} part, suppose that $\gamma$ is connected with
$\delta_0$ by an edge of $\tts$ for an encoding $(\gamma,\varepsilon)$ of $\gamma_0$.
For every such an encoding the vertex $\gamma$ has the form 
$\gamma{\eeq}\{\vv{C_0},\vv{C_1}\}$ for some non-separating circle $C_1$,
and since $\gamma$ is connected by an edge with $\delta_0$, both circles $C_0$, $C_1$
are isotopic to circles disjoint from $D_0$. 
We may assume that the circles $C_0$, $C_1$ themselves are disjoint from $D_0$. 
Then, by Lemma \ref{same-part}, either both circles $C_0$, $C_1$ are contained in $S_0$, 
or both of them are contained in $S_1$. 
Since the subsurface $S_0$ is a torus with one hole, 
it does not contain any bounding pairs of circles. 
Since $C_0$, $C_1$ is a bounding pair of circles,
this implies that $C_0$, $C_1$ are contained in $S_1$. 
In particular, $C_0$ is a circle in $S_1$. 
This proves the \emph{``if''} direction of the lemma, and hence completes the proof. \eproof

\mytitle{Remark.} The assumption $\genus(S)\geqs 4$ in Lemmas \ref{c-in-t-encode} 
and \ref{c-out-t-encode} is needed only to ensure the existence of admissible encodings. 
If $\genus(S)\leq 3$, then there are no marked triangles in $\tts$ 
and hence there are no admissible encodings.

\mypar{Encodings and edges of $\ccs$.}{encoding-edges} The goal of this subsection
is to find properties of encodings of two arbitrary vertices of $\ccs$
which are equivalent to the property of being connected by an edge of $\ccs$
and can be stated entirely in terms of the Torelli building $\tts$.
There are several cases to consider.

If both vertices are the isotopy classes of separating circles,
then their encodings are these circles themselves considered as vertices,
and by the definition of $\tts$ they are connected by an edge in $\tts$
if and only if they are connected by an edge in $\ccs$.

\mylemma{Lemma.}{c-d1-edges} \emph{Suppose that $\genus(S)\geqs 4$. 
Let\/ $C_0$, $D_0$ be, respectively, a non-separating and a separating 
circle in\/ $S$, and let\/ $\gamma_0{\eeq}\vv{C_0}$, $\delta_0{\eeq}\vv{D_0}$.
Suppose that\/ $D_0$ is a genus\/ $1$ circle in\/ $S$.}
 
\vspace{-1ex}
\emph{Then\/ $\gamma_0$ is connected with\/ $\delta_0$ by an edge of\/ $\ccs$
if and only if either the condition\/ $({\mathcal S_0})$ of Lemma \ref{c-in-t-encode},
or the condition\/ $({\mathcal S_1})$ of Lemma \ref{c-out-t-encode} holds.}

\proof This is an easy corollary of Lemmas \ref{c-in-t-encode} and \ref{c-out-t-encode}.
Since $D_0$ is a genus $1$ circle, we are in the situation of the subsection \ref{genus-one-non-sep}.
Let $S_0$ and $S_1$ be the subsurfaces of $S$ defined at the beginnig of the subsection \ref{genus-one-non-sep}.
In particular, $S_0$ is a torus with one hole.

If $\gamma_0$ is connected with $\delta_0$ by an edge of $\ccs$, 
then $C_0$ is isotopic to a circle disjoint from $D_0$, 
and hence to a circle contained in $S_0$ or $S_1$. 
If $C_0$ is isotopic to a circle contained in $S_0$,
Lemma \ref{c-in-t-encode} applies and by this lemma the the condition $({\mathcal S_0})$ holds.
If $C_0$ is isotopic to a circle contained in $S_1$,
Lemma \ref{c-out-t-encode} applies and by this lemma the the condition $({\mathcal S_1})$ holds.
This proves the \emph{``only if''} direction of the lemma.

If the condition $({\mathcal S_0})$ from Lemma \ref{c-in-t-encode} holds,
then Lemma \ref{c-in-t-encode} implies that $C_0$ is isotopic to a circle contained in $S_0$,
and hence $\gamma_0$ is connected with $\delta_0$ by and edge of $\ccs$. 
The argument for $({\mathcal S_1})$ is similar. 
This proves the \emph{``if''} direction of the lemma, 
and hence completes the proof. \eproof

\mylemma{Lemma.}{c-d-edges} \emph{Suppose that $\genus(S)\geqs 4$. 
Let\/ $C_0$, $D_0$ be, respectively, a non-separating and a separating 
circle in\/ $S$, and let\/ $\gamma_0{\eeq}\vv{C_0}$, $\delta_0{\eeq}\vv{D_0}$.
Suppose that\/ $D_0$ is not a genus\/ $1$ circle.}

\vspace{-1ex} 
\emph{Then\/ $\gamma_0$ is connected with\/ $\delta_0$ by an edge of\/ $\ccs$
if and only if there exists an admissible encoding $(\gamma,\varepsilon)$ of $\gamma_0$ 
such that $\gamma$ is connected with $\delta_0$ by an edge of $\tts$.}

\proof Let $S_0$, $S_1$ be two subsurfaces bounded by $D_0$ in $S$. 
Since $D_0$ is not a genus $1$ circle, both $S_0$ and $S_1$ have genus $\geqs 2$. 

If $\gamma_0$ and $\delta_0$ are connected by an edge of $\ccs$, 
then $C_0$ is isotopic to a circle disjoint from $D_0$ and we may assume that
$C_0$ itself is disjoint from $D_0$. Then $C_0$ is contained in one of the subsurfaces
$S_0$, $S_1$. We may assume that $C_0$ is contained in $S_0$.
Since $\genus(S_0)\geqs 2$, there is a non-separating circle $C_1$ in $S_0$
such that $C_0$, $C_1$ is a bounding pair of circles in $S_0$.
Let $\gamma{\eeq}\{\vv{C_0},\vv{C_1}\}$.  
Since the circles $C_0$, $C_1$ are disjoint from $D_0$, 
$\gamma$ is connected with $\delta_0$ by and edge of $\tts$. 
In addition, since $\genus(S)\geq 4$, there exists a vertex $\varepsilon$
such that $(\gamma,\varepsilon)$ is an admissible encoding of $\gamma_0$ 
(it does not matter if $\varepsilon$ is the isotopy class of a circle contained in $S_0$ or not). 
This proves the \emph{``only if''} part of the lemma.

Let us prove the \emph{``if''} part. If $\gamma$ and $\delta_0$ are connected by and edge,
then $\gamma{\eeq}\{\vv{E_0},\vv{E_1}\}$ for some bounding pair of circles $E_0$, $E_1$ disjoint from $D_0$.
Since $(\gamma,\delta)$ encodes $\gamma_0$, the circle $C_0$ is isotopic to one of these circles $E_0$, $E_1$.
In particular, $C_0$ is isotopic to a circle disjoint form $D_0$, i.e.
$\gamma_0$ and $\delta_0$ are connected by an edge of $\ccs$. 
This proves the \emph{``if''} part of the lemma and hence completes the proof. \eproof

\mylemma{Lemma.}{non-sep-non-sep} \emph{Suppose that $\gen{\eeq}\genus(S)\geqs 5$.
Let\/ $\gamma'_0$, $\gamma''_0$ be two isotopy classes of non-separating
circles in\/ $S$. The following two conditions are equivalent.}
\vspace{-2ex}
\begin{itemize}
\item[($ \ccc$)] \emph{$\gamma'_0$ and\/ $\gamma''_0$ are connected by an edge in\/ $\ccs$\/.}
\vspace{-1ex}
\item[($ \ttg$)] \emph{There exist encodings\/ $(\gamma',\delta')$ of\/ $\gamma'_0$ and\/ $(\gamma'',\delta'')$ of\/ $\gamma''_0$
such that\/ $\gamma'$ and\/ $\gamma''$ are connected by an edge in\/ $\tts$.}
\end{itemize}
\vspace{-2ex}

\emph{Proof.} Let $\gamma'_0{\eeq}\vv{C'_0}$ and $\gamma''_0{\eeq}\vv{C''_0}$ for some non-separating circles
$C'_0$ and $C''_0$ in $S$.

\textbf{($ \ccc$) \textit{implies}\/ ($ \ttg$).}\/ Suppose first that the condition ($ \ccc$) holds. Then we may assume that the circles $
C'_0$ and $C''_0$ are disjoint and are not isotopic (the case of isotopic circles being trivial). 
There are two cases to consider.

\mytitle{Case 1. $C'_0\cup C''_0$ separates $S$.} Then $C'_0$, $C''_0$ is a bounding pair of
circles. Let $\gamma{\eeq}\{\vv{C'_0},\vv{C''_0}\}$ be the corresponding vertex of $\tts$. 
Then for some separating vertices $\delta'$, $\delta''$ of $\tts$ 
the pairs $(\gamma,\delta')$ and $(\gamma,\delta'')$ are encoding $\gamma'_0$ and $\gamma''_0$
respectively. This proves that ($ \ccc$) implies ($ \ttg$) in Case 1. \esubproof

\mytitle{Case 2. $C'_0\cup C''_0$ does not separate $S$.} Suppose that the genus $\gen$ of $S$ is fixed.
Then up to a diffeomorphism there is only one (ordered) triple $(S, C'_0, C''_0)$ such that the union $C'_0\cup C''_0$ n
ot separating $S$ (and such a triple exists if and only if $\gen\geqs 2$). 
Therefore, we can deal with this case by considering a model example.

Since $\gen\geqs 5$ by the assumption, $\gen{\eeq}\gen'+\gen''+2$ for some $\gen'\geqs 1$ and $\gen''\geqs 2$.
Let $R'$ be a surface of genus $\gen'{\eeq}\genus(R')\geqs 1$ having exactly two boundary components $C'_0$, $C'_1$, and
let $R''$ be a surface of genus $\gen''{\eeq}\genus(R'')\geqs 2$ having exactly two boundary components $C''_0$, $C''_1$.
Let $R_0$ be a surface of genus $\genus(R_0){\eeq}0$ having exactly $4$ boundary componets. 
Let us glue $R'$ and $R''$ to $R_0$ by glueing $\partial R'$ to arbitrary $2$ of the $4$ boundary components of $R_0$,
and glueing $\partial R''$ to the $2$ other boundary components of $R_0$.
The resulting surface $R$ is a closed surface of genus $\gen'+\gen''+2$.
We may take $(R, C'_0, C''_0)$ as our model of $(S, C'_0, C''_0)$.

Clearly, $\{C'_0,C'_1\hspace{0.05em}\}$ and $\{C''_0,C''_1\hspace{0.05em}\}$ are bounding pairs of circles.
Let $\gamma'$ and $\gamma''$ be the corresponding vertices of $\tts$.
Since the circles $C'_0$, $C'_1$, $C''_0$, $C''_1$ are pair-wise disjoint by the construction,
the vertices $\gamma'$ and $\gamma''$ are connected by an edge of\/ $\tts$.
On the other hand, for some vertices $\delta'$, $\delta''$ of\/ $\tts$ corresponding to separating circles
the pairs $(\gamma,\delta')$ and $(\gamma,\delta'')$ are encoding $\gamma'_0$ and $\gamma''_0$
respectively. This proves that ($ \ccc$) implies ($ \ttg$) in Case 2. \esubproof

\textbf{($ \ttg$) \textit{implies}\/ ($ \ccc$).}\/ Suppose now that the condition ($ \ttg$) holds.

Let $(\gamma',\delta')$ and $(\gamma'',\delta'')$ 
be encodings of $\gamma'_0$ and $\gamma''_0$ respectively
such that $\gamma'$ and $\gamma''$ are connected by an edge in $\tts$.
Then $\gamma'{\eeq}\{\vv{C'_0},\vv{C'_1}\}$ and $\gamma''{\eeq}\{\vv{C''_0},\vv{C''_1}\}$
for some bounding pairs of circles $\{C''_0,C''_1\hspace{0.05em}\}$ and $\{C''_0,C''_1\hspace{0.05em}\}$.
Since $\gamma'$ is connected by an edge of $\tts$ with $\gamma''$,
we may assume that circles $C'_0$, $C'_1$, $C''_0$, $C''_1$ are pair-wise disjoint.
Then, in particular, the circles $C'_0$ and $C''_0$ are disjoint, 
and hence $\gamma'0$ and $\gamma''_0$ are connected by an edge in $\ccs$.
This proves that ($ \ttg$) implies ($ \ccc$), and hence completes the proof of the lemma. \eproof

\mypar{Proof of Theorem \ref{automorphisms of geometry}.}{proof-aut-tg}
Let $S$ be a closed surface of genus $\geqs\, 5$. 
Let $A\colon\tts\tto\tts$ be an automorphism of the Torelli
geometry $\tts$. 

Recall that the set of vertices of $\ccs$ which are the isotopy classes of separating circles
is equal to the set of $\mathcal{SC}$-vertices of $\tts$, 
i.e. the set of vertices of $\tts$ marked by the symbol $\mathcal{SC}$.
Since $A$ is an automorphism of the Torelli building $\tts$,\sss $A$ preserves markings of vertices.
Therefore, $A$ induces a bijective self-map of
the set of vertices of $\ccs$ which are the isotopy classes of separating circles.
  
By Corollaries \ref{pairs-are-enc} and \ref{encoding}, encodings by admissible pairs establish
a canonical $1$-to-$1$ correspondence between the isotopy classes of non-separating circles
and the equivalence classes of admissible pairs with respect to
the relation of being connected by a sequence of moves (obviously, this is an equivalence relation).
Since the admissible pairs and the moves of both types are defined entirely in terms of the Torelli building,
$A$ maps admissible pairs to admissible pairs, and if two admissible pairs are related by a move, 
their images under $A$ are also related by a move of the same type.
Therefore, $A$ induces bijective self-maps of the set of admissible pairs respecting
the relation of being connected by a sequence of moves.
It follows that $A$ induces also a bijective self-map of the set of the equivalence classes 
of admissible pairs with respect to the above relation.
When combined with encodings, the last self-map induced by $A$
induces, in turn, a bijective self-map of the set of the isotopy classes of non-separating circles.
In other words, $A$ induces a bijective self-map of the set of vertices of $\ccs$ 
which are the isotopy classes of non-separating circles.

By combining the results of the last paragraphs, we see that $A$ induces
a bijective self-map of the set of the isotopy classes of non-trivial circles in $S$,
i.e. of the set of the vertices of $\ccs$.
We will denote this induced self-map by $A_\mathcal{C}$. 
By the construction, the induced map $A_\mathcal{C}$ agrees with $A$ and encodings 
of non-separating vertices the sense of the following lemma.  

\mylemma{Lemma}{a-encod} \emph{If $(\gamma,\delta)$ is an admissible encoding 
of a non-separating vertex $\gamma_0$ of $\ccs$,
then $A_\mathcal{C}(\gamma_0)$ is also a non-separating vertex of $\ccs$,
and $(A(\gamma),A(\delta))$ is an admissible encoding of $A_\mathcal{C}(\gamma_0)$.} \eproof.

Now we have to prove that $A_\mathcal{C}$ is not only a bijective self-map of the set of vertices
of $\ccs$, but, moreover, is an automorphism of the complex of curves $\ccs$. 
This is done in the following lemma.
 
\mylemma{Lemma.}{edges-preserved} \emph{$A_\mathcal{C}$ preserves the structure of the simplicial complex $\ccs$.
In other words, $A_\mathcal{C}$ maps every set of vertices which is a simplex of $\ccs$ to a set 
which is also a simplex of $\ccs$.}

\proof Since $\ccs$ is a flag complex (see the subsection \ref{cc}), 
it is sufficient to prove that if two vertices of $\ccs$ are connected by an edge,
then their images under the induced map $A_\mathcal{C}$ are also connected by an edge.
There a several cases to consider.

\mytitle{Case 1. Two separating vertices.} By the definitions of $\ccs$ and $\tts$, two separating vertices, 
i.e. two isotopy classes of non-trivial separating circles, 
are connected by an edge in $\ccs$ if and only if they are connected by an edge in $\tts$.
It follows that $A_\mathcal{C}$ takes edges connecting two such vertices into similar edges. \esubproof

\mytitle{Case 2. Two non-separating vertices.} Let $\gamma'_0$, $\gamma''_0$ be two non-separating vertices
of $\ccs$. By Lemma \ref{non-sep-non-sep} the vertices $\gamma'_0$ and $\gamma''_0$ are connected by an edge of $\ccs$
if and only the condition ($ \ttg$) of this lemma holds. Lemma \ref{a-encod} implies that if this condition
holds for $\gamma'_0$, $\gamma''_0$, then it holds also for $A_\mathcal{C}(\gamma'_0)$, $A_\mathcal{C}(\gamma''_0)$
in the roles of $\gamma'_0$, $\gamma''_0$ respecively. By using Lemma \ref{non-sep-non-sep} again,
we conclude that if $\gamma'_0$, $\gamma''_0$ are connected by and edge, 
then $A_\mathcal{C}(\gamma'_0)$, $A_\mathcal{C}(\gamma''_0)$ are also connected by an edge. \esubproof

\mytitle{Case 3. A non-separating vertex and a separating vertex.} Let\/ $C_0$, $D_0$ be, 
respectively, a non-separating and a separating 
circle in\/ $S$, and let\/ $\gamma_0{\eeq}\vv{C_0}$, $\delta_0{\eeq}\vv{D_0}$.
There are two subcases depending on if $D_0$ is a genus $1$ circle or not.

First of all, by Lemma \ref{links} $D_0$ is a genus $1$ circle if and only if 
the graph $\link\sep(\delta_0)$ is connected. In view of the definition of $\link\sep(\delta_0)$ 
and Lemma \ref{a-encod}, if $\link\sep(\delta_0)$ is connected, 
then $\link\sep(A_\mathcal{C}(\delta))$ is also connected.
It follows that $A_\mathcal{C}$ preserves the set of separating vertices of $\ccs$ 
which are the isotopy classes of genus $1$ circles.

\mytitle{Subcase 3-A. $D_0$ is a genus $1$ circle.}  As we just saw,
then $A_\mathcal{C}(\delta)$ is also an isotopy class of genus $1$ circle.
By Lemma \ref{c-d1-edges}, $\gamma_0$ is connected with $\delta_0$ by an edge of $\ccs$ 
if and only if either the condition\/ $({\mathcal S_0})$ of Lemma \ref{c-in-t-encode},
or the condition\/ $({\mathcal S_1})$ of Lemma \ref{c-out-t-encode} holds.
By arguments completely similar to the arguments in Cases 1 and 2, 
we conclude that if $\gamma_0$, $\delta_0$ are connected by an edge,
then $A_\mathcal{C}(\gamma_0)$, $A_\mathcal{C}(\delta_0)$ are also connected by an edge. \esubproof 

\mytitle{Subcase 3-B. $D_0$ is not a genus $1$ circle.} This subcase is
similar to the previous one. The only difference is that the 
reference to Lemma \ref{c-d1-edges} should be replaced by a reference to Lemma \ref{c-d-edges}. 
This completes the proof of the lemma. \eproof

\mytitle{End of the proof of Theorem \ref{automorphisms of geometry}.}
By the last lemma  $A_\mathcal{C}$ is an automorphism of $\ccs$.
Therefore, Theorem \ref{aut-cc} implies that $A_\mathcal{C}{\eeq}\ffc$ 
for some diffeomorphism $F\colon S\rightarrow S$.
Recall that for every non-trivial circle $C$, the map $\ffc$ maps
the vertex $\llvv{C}$ of $\ccs$ to the vertex $\aavv{F(C)}$,
and that $\ffc$ is a bijective self-map of the set of all vertices of $\ccs$.
See the subsection \ref{cc}.
The bijection $\ffc$, in turn, induces an automorphism $\fft\colon\tts{\tto}\tts$. 
See the subsection \ref{torelli-geom} for the details.
In order to complete the proof of the theorem, it is sufficient to prove that $A{\eeq}\fft$. 

Let $\gamma{\eeq} \{\gamma_0,\gamma_1\}$ be an arbitrary $\mathcal{BP}$-vertex of $\tts$.
Then there exists an admissible pair $(\gamma, \delta_0)$ encoding the vertex $\gamma_0$ of $\ccs$.
By Lemma \ref{a-encod}, this implies $(A(\gamma), A(\delta_0))$ is an encoding of $A_\mathcal{C}(\gamma_0)$.
But $A_\mathcal{C}{\eeq}\ffc$ and hence $A_\mathcal{C}(\gamma_0){\eeq}\ffc (\gamma_0)$.
It follows that $(A(\gamma), A(\delta_0))$ is an encoding of $\ffc (\gamma_0)$.
In particular, $\ffc (\gamma_0)$ is one of the elements of the bounding pair $A(\gamma)$
(recall that bounding pairs are, in particular, $2$-element sets of vertices of $\ccs$).
By a similar argument, $\ffc (\gamma_1)$ is also one of the elements of the bounding pair $A(\gamma)$.
But the isotopy classes $\ffc (\gamma_0)$, $\ffc (\gamma_1)$ are different, 
because the isotopy classe $\gamma_0$, $\gamma_1$ are different, being elements of a bounding pair
(namely, of the bounding pair $\gamma$). Therefore, the pair $A(\gamma)$ is equal 
to the pair $\{\hspace{0.05em}\ffc (\gamma_0),\,\ffc (\gamma_1)\hspace{0.05em}\}$.
The latter pair is nothing else as $\fft (\gamma)$.
It follows that $A(\gamma){\eeq}\fft (\gamma)$.

$C_0$ is not isotopic to $C_1$, and hence $F(C_0)$ is not isotopic to $F(C_1)$.
It follows that the pair $A(\gamma)$ should be equal to the pair $\{\hspace{0.05em}\ffc (\gamma_0),\,\ffc (\gamma_1)\hspace{0.05em}\}$.
The latter pair is nothing else as $\fft (\gamma)$.
It follows that $A(\gamma){\eeq}\fft (\gamma)$.

Since $\gamma$ was an arbitrary $\mathcal{BP}$-vertex of $\tts$, this means $A$ agrees with $\fft$
on the set of $\mathcal{BP}$-vertices. Since, as we noted above, $A$ agrees with $\fft$ also
on the set of $\mathcal{SC}$-vertices, this implies that $A{\eeq}\fft$, and hence completes the proof 
of the theorem. \eproof

\begin{flushright}

October 9, 2014

\vspace{\bigskipamount}

Benson Farb:\hspace*{3.25em}
Dept. of Mathematics, University of Chicago{ }\\
5734 University Ave.\hspace*{10.7em}{ }\\
Chicago, Il 60637\hspace*{12.24em}{ }\\
E-mail: farb@math.uchicago.edu\hspace*{5.15em}{ }

\vspace{\bigskipamount}

Nikolai V. Ivanov:\hspace*{1em} 
http:/\!/\!nikolaivivanov.com\hspace*{8.25em}{ }\\
E-mail: nikolai@nikolaivivanov.com\hspace{4.0em}

\end{flushright}

\end{document}